\documentclass[a4paper,11pt]{article}
\usepackage{amsmath, amssymb, amsthm}
\usepackage{graphicx} 

\usepackage[latin1]{inputenc}

\begin{document}

\def\prod{{\rm prod}} \def\R{{\mathbb R}} \def\cos{{\rm cos}}
\def\lg{{\rm lg}} \def\SS{{\mathcal S}} \def\PP{{\mathcal P}}
\def\CC{{\mathcal C}} \def\Cox{{\rm Cox}} \def\Art{{\rm Art}}
\def\AA{{\mathcal A}} \def\UU{{\mathcal U}} \def\Link{{\rm Link}}
\def\E{{\mathbb E}} \def\LL{{\mathcal L}} \def\VV{{\mathcal V}}
\def\cst{{\rm cst}} \def\Star{{\rm Star}} \def\Span{{\rm Span}}
\def\dim{{\rm dim}} \def\End{{\rm End}} \def\SSS{{\mathfrak S}}
\def\TT{{\mathcal T}} \def\link{{\rm link}}


\title{\bf{$K(\pi,1)$ and word problems for infinite type Artin-Tits groups, and applications to virtual braid groups}}

\author{\textsc{Eddy Godelle\footnotemark[1] and Luis Paris\footnote{Both authors are partially supported by the {\it Agence Nationale de la Recherche} ({\it projet Théorie de Garside}, ANR-08-BLAN-0269-03).}}}

\date{\today}

\maketitle

\begin{abstract}
\noindent
Let $\Gamma$ be a Coxeter graph, let $(W,S)$ be its associated Coxeter system, and let $(A,\Sigma$) be its associated Artin-Tits system. We regard $W$ as a reflection group acting on a real vector space $V$. Let $I$ be the Tits cone, and let $E_\Gamma$ be the complement in $I +iV$ of the reflecting hyperplanes. Recall that Charney, Davis, and Salvetti have constructed a simplicial complex $\Omega(\Gamma)$ having the same homotopy type as $E_\Gamma$. We observe that, if $T \subset S$, then $\Omega(\Gamma_T)$ naturally embeds into $\Omega (\Gamma)$. We prove that this embedding admits a retraction $\pi_T: \Omega(\Gamma) \to \Omega (\Gamma_T)$, and we deduce several topological and combinatorial results on parabolic subgroups of $A$. From a family $\SS$ of subsets of $S$ having certain properties, we construct a cube complex $\Phi$, we show that $\Phi$ has the same homotopy type as the universal cover of $E_\Gamma$, and we prove that $\Phi$ is  CAT(0) if and only if $\SS$ is a flag complex. We say that $X \subset S$ is free of infinity if $\Gamma_X$ has no edge labeled by $\infty$. We show that, if $E_{\Gamma_X}$ is aspherical and $A_X$ has a solution to the word problem for all $X \subset S$ free of infinity, then $E_\Gamma$ is aspherical and $A$ has a solution to the word problem. We apply these results to the virtual braid group $VB_n$. In particular, we give a solution to the word problem in $VB_n$, and we prove that the virtual cohomological dimension of $VB_n$ is $n-1$. 
\end{abstract}





\section{Introduction}

We start with some basic definitions on Coxeter groups and Artin-Tits groups. Let $S$ be a finite set. A {\it Coxeter matrix} over $S$ is a square matrix $M=(m_{s,t})_{s,t \in S}$ indexed by the elements of $S$, such that $m_{s,s}=1$ for all $s \in S$, and $m_{s,t} = m_{t,s} \in \{2,3,4, \dots, \infty\}$ for all $s,t \in S$, $s \neq t$. Such a Coxeter matrix is usually represented by its {\it Coxeter graph}, denoted $\Gamma$. This is a labeled graph whose set of vertices is $S$, where two vertices $s$ and $t$ are joined by an edge if $m_{s,t} \ge 3$, and where this edge is labeled by $m_{s,t}$ if $m_{s,t} \ge 4$.

\bigskip\noindent
The {\it Coxeter system} of $\Gamma$ is the pair $(W,S)= (W_\Gamma, S)$, where $W$ is the group 
\[
W= \left\langle S\ \left|\
\begin{array}{ll}
s^2=1 &\text{ for all } s \in S\\
(st)^{m_{s,t}}=1 &\text{ for all } s,t \in S,\ s \neq t,\text{ and } m_{s,t} \neq \infty
\end{array}\right.\right\rangle\,.
\]
The group $W$ is called the {\it Coxeter group} of $\Gamma$.

\bigskip\noindent 
If $a,b$ are two letters and $m$ is an integer greater or equal to $2$, we set $\prod(a,b:m) =(ab)^{\frac{m}{2}}$ if $m$ is even, and $\prod(a,b:m) =(ab)^{\frac{m-1}{2}}a$ if $m$ is odd. We take $\Sigma = \{ \sigma_s; s \in S\}$, a set in one-to-one correspondence with $S$. The {\it Artin-Tits system} of $\Gamma$ is the pair $(A,\Sigma) = (A_\Gamma, \Sigma)$, where $A$ is the group
\[
A= \langle \Sigma\ |\ \prod(\sigma_s, \sigma_t : m_{s,t}) = \prod(\sigma_t, \sigma_s : m_{s,t}) \text{ for } s,t \in S,\ s \neq t, \text{ and } m_{s,t} \neq \infty \rangle\,.
\]
The group $A$ is called the {\it Artin-Tits group} of $\Gamma$.

\bigskip\noindent
The map $\Sigma \to S$ which sends $\sigma_s$ to $s$ for all $s \in S$ induces an epimorphism $\theta : A \to W$. The kernel of $\theta$ is called the {\it colored Artin-Tits group} of $\Gamma$ and is denoted by $CA=CA_\Gamma$. On the other hand, $\theta : A \to W$ has a ``natural'' set-section $\tau: W \to A$ defined as follows. Let $w \in W$. Let $w=s_1 s_2 \cdots s_l$ be a reduced expression for $w$. Then $\tau(w) = \sigma_{s_1} \sigma_{s_2} \cdots \sigma_{s_l}$. Tits' solution to the word problem for Coxeter groups (see \cite{Tits1}) implies that the definition of $\tau(w)$ does not depend on the choice of the reduced expression. 

\bigskip\noindent
The Coxeter group $W$ has a faithful linear representation $W \hookrightarrow GL(V)$,  called {\it canonical representation}, where $V$ is a real vector space of dimension $|S|$. The group $W$, viewed as a subgroup of $GL(V)$, is generated by reflections, and acts properly discontinuously on an convex open cone $I$, called {\it Tits cone} (see \cite{Bourb1}). The set of reflections in $W$ is $R= \{w s w^{-1}; s \in S \text{ and } w \in W\}$, and $W$ acts freely and properly discontinuously on $I \setminus (\cup_{r \in R} H_r)$, where, for $r \in R$, $H_r$ denotes the hyperplane of $V$ fixed by $r$.  Set
\[
E= E_\Gamma = (I \times V) \setminus \left( \bigcup_{r \in R} H_r \times H_r \right)\,.
\]
This is a connected manifold of dimension $2|S|$ on which the group $W$ acts freely and properly discontinuously. One of the main result in the subject is the following.

\bigskip\noindent
{\bf Theorem} (Van der Lek \cite{Lek1}). 
{\it $\pi_1(E_\Gamma) = CA_\Gamma$, $\pi_1(E_\Gamma/W) = A_\Gamma$, and the exact sequence associated to the regular cover $E_\Gamma \to E_\Gamma/W$ is $1 \to CA_\Gamma \to A_\Gamma \to W \to 1$.}

\bigskip\noindent
We say that $\Gamma$ is of {\it type $K(\pi,1)$} if $E=E_\Gamma$ is an Eilenberg MacLane space, that is, if $E=E_\Gamma$ is a $K(CA_\Gamma,1)$ space. A central conjecture in the theory of Artin-Tits groups, known as the {\it $K(\pi,1)$ conjecture}, is that all the Coxeter graphs are of type $K(\pi,1)$. 

\bigskip\noindent
Artin-Tits groups are badly understood in general. In particular, it is not known whether they all have solutions to the word problem, and it is not know whether they all are torsion free. Actually, the theory of Artin-Tits groups mainly consists on the study of some more or less extended families of Artin-Tits groups. 

\bigskip\noindent
The first family of Artin-Tits groups which has been studied is the family of spherical type Artin-Tits groups (see \cite{Bries1}, \cite{Bries2}, \cite{BriSai1}, \cite{Delig1}). Recall that an Artin-Tits group $A_\Gamma$ is of {\it spherical type} if its associated Coxeter group $W_\Gamma$ is finite. Spherical type Artin-Tits groups are torsion free,  they have (fast) solutions to the word problem (see \cite{BriSai1}, \cite{Delig1}, \cite{Charn1}, \cite{Charn2}), and they are of type $K(\pi,1)$ (see \cite{Delig1}). They are also of importance in the study of other sorts of Artin-Tits groups.

\bigskip\noindent
For $X\subset S$, set $M_X = (m_{s,t})_{s,t \in X}$, denote by $\Gamma_X$ the Coxeter graph of $M_X$, by $W_X$ the subgroup of $W=W_\Gamma$ generated by $X$, set $\Sigma_X = \{ \sigma_s; s \in X\}$, and denote by $A_X$ the subgroup of $A = A_\Gamma$ generated by $\Sigma_X$. By \cite{Bourb1}, the pair $(W_X,X)$ is a Coxeter system of $\Gamma_X$, and, by \cite{Lek1}, the pair $(A_X,\Sigma_X)$ is an Artin-Tits system of $\Gamma_X$. A new proof of this last result will be given in Section~2. The subgroup $W_X$ is called {\it standard parabolic subgroup} of $W$, and $A_X$ is called {\it standard parabolic subgroup} of $A$.

\bigskip\noindent
Two families of subsets of $S$ will play a major role in the paper. The first family, denoted by $\SS^f$, is made of the subsets $X \subset S$ such that $W_X$ is finite. For $X \subset S$, we say that $\Gamma_X$ is {\it free of infinity} if $m_{s,t} \neq \infty$ for all $s,t \in X$. The second family, denoted by $\SS_{< \infty}$, is made of the subsets $X \subset S$ such that $\Gamma_X$ is free of infinity. Note that $\SS^f \subset \SS_{<\infty}$.

\bigskip\noindent
We say that $\Gamma$ (or $A_\Gamma$) is of {\it FC type} if $\SS^f = \SS_{<\infty}$. In other words, $\Gamma$ is of FC type if $\SS^f$, viewed as an abstract simplicial complex over $S$, is a flag complex (note that FC stands for ``flag complex''). FC type Artin-Tits groups form another family of extensively studied Artin-Tits groups. One of the main results concerning this family is about a cube complex $\Phi= \Phi(\Gamma, \SS^f)$ constructed by Charney and Davis in \cite{ChaDav1}: it is proved in \cite{ChaDav1} that $\Phi$ has always the same homotopy type as the universal cover of $E_\Gamma$, and that $\Phi$ endowed with its cubical structure is CAT(0) if and only if $\Gamma$ is of FC type. It is well-known that CAT(0) metric spaces are contractible (see \cite{BriHae1}), thus this implies that FC type Artin-Tits groups are of type $K(\pi,1)$. The cubical structure on $\Phi$ is also used by Altobelli and Charney \cite{AltCha1} to solve the word problem in FC type Artin-Tits groups.

\bigskip\noindent
We say that a family $\SS$ of subsets of $S$ is {\it complete and $K(\pi,1)$} if (1) $\SS$ is closed under inclusion (that is, if $X \in \SS$ and $Y \subset X$, then $Y \in \SS$), (2) $\Gamma_X$ is of type $K(\pi,1)$ for all $X \in \SS$, (3) $\SS^f \subset \SS$. The starting idea in the present paper consists on replacing the family $\SS^f$ in the study of Charney and Davis by a complete and $K(\pi,1)$ family. To such a family $\SS$ we associate a cube complex $\Phi = \Phi(\Gamma, \SS)$, we prove that $\Phi$ has always the same homotopy type as the universal cover of $E_\Gamma$, and we show that $\Phi$ is CAT(0) if and only if $\SS$, viewed as an abstract simplicial complex over $S$, is a flag complex.

\bigskip\noindent
If $\Gamma_X$ is of type $K(\pi,1)$ for all $X \in \SS_{<\infty}$, then $\SS_{<\infty}$ is complete and $K(\pi,1)$ and is a flag complex. Reciprocally, if $\SS$ is a complete and $K(\pi,1)$ family and is a flag complex, then $\SS_{<\infty} \subset \SS$ and $\Gamma_X$ is of type $K(\pi,1)$ for all $X \in \SS_{<\infty}$ (this is explained in more details in Section 4). So, the study of the families of subsets of $S$ that are complete and $K(\pi,1)$ and that are flag complexes can be restricted without lost of generality to the study of $\SS_{<\infty}$ under the assumption that $\Gamma_X$ is of type $K(\pi,1)$ for all $X \in \SS_{<\infty}$.

\bigskip\noindent
In Section 5 we assume that $\Gamma_X$ is of type $K(\pi,1)$ and $A_X$ has a solution to the word problem for all $X \in \SS_{<\infty}$, and we use the geometry of $\Phi = \Phi(\Gamma, \SS_{<\infty})$ to solve the word problem in $A$.

\bigskip\noindent
In the last section we apply the previous results to the virtual braid group $VB_n$. It is known that $VB_n$ is a semi-direct product $VB_n = K_n \rtimes \SSS_n$ of an Artin-Tits group $K_n$ with the symmetric group $\SSS_n$ (see \cite{Raben1}, \cite{BarBel1}). Let $\Gamma_{VB,n}$ denote the Coxeter graph of $K_n$ and let $\SS_{<\infty}$ be the set of subsets $X \subset S$ such that $(\Gamma_{VB,n})_X$ is free of infinity. We prove that $(\Gamma_{VB,n})_X$ is of type $K(\pi,1)$ and $(K_n)_X$ has a solution to the word problem for all $X \in \SS_{<\infty}$. We deduce that the cohomological dimension of $K_n$ is $n-1$, the virtual cohomological dimension of $VB_n$ is $n-1$, and $VB_n$ has a solution to the word problem. 

\bigskip\noindent
Although our work is inspired by \cite{ChaDav1} and \cite{AltCha1}, most of our proofs are different. For some results we prefer to provide new proofs that we think are simpler. However, for other results, the proofs in \cite{ChaDav1} and/or in \cite{AltCha1} cannot be extended because  they are based on results on spherical type Artin-Tits groups that are wrong for the other Artin-Tits groups. To overcome this problem we develop new tools in Section 2.

\bigskip\noindent
In Section 2 we start with a simplicial complex $\Omega = \Omega(\Gamma)$, introduced by Charney and Davis \cite{ChaDav2} and, independently, by Salvetti \cite{Salve1}, which has the same homotopy type as $E=E_\Gamma$. If $T \subset S$, then $\Omega({\Gamma_T})$ embeds in a natural way into $\Omega(\Gamma)$. The main tool of this paper, proved at the beginning of Section 2, is the following.

\bigskip\noindent 
{\bf Theorem 2.2.}
{\it The natural embedding $\Omega (\Gamma_T) \hookrightarrow \Omega(\Gamma)$ admits a retraction $\pi_T : \Omega(\Gamma) \to \Omega({\Gamma_T})$.} 

\bigskip\noindent
The study of this retraction gives rise to some results on standard parabolic subgroups of $A_\Gamma$. Firstly, a straightforward consequence of the existence of this retraction is the following result, which, curiously, was unknown before: if $T \subset S$ and $\Gamma$ is of type $K(\pi,1)$, then $\Gamma_T$ is also of type $K(\pi,1)$ (see Corollary 2.4). We also give new proofs of some known results by Van der Lek \cite{Lek1}: for $T \subset S$, $(A_T, \Sigma_T)$ is an Artin-Tits system of $\Gamma_T$, and, for $R,T \subset S$, we have $A_R \cap A_T = A_{R \cap T}$ (see Theorem 2.5). Finally, under the assumption that $A$ has a solution to the word problem, we give an algorithm which, given $T \subset S$, $\alpha \in A$, and a word $\omega \in (\Sigma \sqcup \Sigma^{-1})^\ast$ which represents $\alpha$, decides whether $\alpha \in A_T$, and, if yes, determines a word $\kappa_T(\omega) \in (\Sigma_T \sqcup \Sigma_T^{-1})^\ast$ which represents $\alpha$ (see Proposition 2.7).

\bigskip\noindent
The paper is organized as follows. In Section 2 we study the retraction $\pi_T : \Omega(\Gamma) \to \Omega(\Gamma_T)$ and its different applications. In Section 3 we define the complex $\Phi= \Phi(\Gamma,\SS)$ and prove that, if $\SS$ is a complete and $K(\pi,1)$ family of subsets of $S$, then $\Phi$ has the same homotopy type as the universal cover of $E_\Gamma$ (Theorem 3.1). In Section 4 we define a cubical structure on $\Phi$ and we show that $\Phi$ is a CAT(0) cube complex if and only if $\SS$ is a flag complex (Theorem 4.2). In Section 5 we prove that the group $A_\Gamma$ has a solution to the word problem if $\Gamma_X$ is of type $K(\pi,1)$ and $A_X$ has a solution to the word problem for all $X \in \SS_{<\infty}$ (Theorem 5.6). In Section~6 we apply the previous results to the virtual braid groups. 




\section{Topology and combinatorics of parabolic subgroups}

We keep the notations of Section 1. So, $\Gamma$ is a Coxeter graph, $(W,S)$ is the Coxeter system of $\Gamma$, and $(A,\Sigma) = (A_\Gamma, \Sigma)$ is the Artin-Tits system of $\Gamma$. We assume $W$ to be embedded into $GL(V)$ via the canonical representation, where $V$ is a real vector space of dimension $|S|$. Let $I$ be the Tits cone (see \cite{Bourb1} for the definition), let $R=\{wsw^{-1}; s \in S \text{ and } w \in W\}$ be the set of reflections in $W$, and let
\[
E=E_\Gamma = (I \times V) \setminus \left( \bigcup_{r \in R} H_r \times H_r \right)\,,
\]
where, for $r \in R$, $H_r$ denotes the hyperplane of $V$ fixed by $r$. Recall that, by \cite{Lek1}, $CA$ is the fundamental group of  $E$ and $A=A_\Gamma$ is the fundamental group of $E/W$. 

\bigskip\noindent
Define the {\it length} of an element $w \in W$ to be the shortest length of a word on $S$ representing $w$. It will be denoted by $\lg (w)$. Let $X,Y$ be two subsets of $S$, and let $w \in W$. We say that $w$ is {\it $(X,Y)$-reduced} if it is of minimal length among the elements of the double-coset $W_X w W_Y$. On the other hand, recall that $\SS^f$ denotes the set of $X \subset S$ such that $W_X$ is finite. We define a relation $\le$ on $W \times \SS^f$ by: $(u,X) \le (v,Y)$ if $X \subset Y$, $v^{-1}u \in W_Y$, and $v^{-1}u$ is  $(\emptyset,X)$-reduced. It is easily checked that $\le$ is an order relation. The chains in $W \times \SS^f$ form an (abstract) simplicial complex called {\it derived complex} of $W \times \SS^f$ and denoted by $(W \times \SS^f)'$. The geometric realization of $(W \times \SS^f)'$ is the { \it Charney-Davis-Salvetti complex}. It will be denoted by $\Omega = \Omega(\Gamma)$. As mentioned before, this complex has been introduced by Charney and Davis in \cite{ChaDav2}, and independently by Salvetti in \cite{Salve1} (Salvetti's construction is for spherical type Artin-Tits groups, but it can be easily extended to the general case). 

\bigskip\noindent
{\bf Theorem 2.1} (Charney, Davis \cite{ChaDav2}, Salvetti \cite{Salve1}).
{\it There exists a homotopy equivalence $\Omega \to  E$ which is equivariant under the action of $W$, and which induces a homotopy equivalence $\Omega/W \to E/W$.}

\bigskip\noindent
Let $T$ be a subset of $S$. Set $\SS^f_T=\{X \in \SS^f; X \subset T\}$. Observe that the inclusion $(W_T \times \SS^f_T) \hookrightarrow (W \times \SS^f)$ induces an embedding $\Omega(\Gamma_T) \hookrightarrow \Omega(\Gamma)$ which is equivariant under the action of $W_T$. The main result of the present section, which will be also the main tool in the remainder of the paper, is the following. 

\bigskip\noindent
{\bf Theorem 2.2.}
{\it \begin{enumerate}
\item
Let $T$ be a subset of $S$. Then the embedding $\Omega(\Gamma_T) \hookrightarrow \Omega(\Gamma)$ admits a retraction  
\break
$\pi_T: \Omega(\Gamma) \to \Omega(\Gamma_T)$ which is equivariant under the action of $W_T$.
\item
Let $R$ and $T$ be two subsets of $S$. Then the restriction of $\pi_T$ to $\Omega(\Gamma_R)$ coincides with $\pi_{R \cap T}: \Omega (\Gamma_R) \to \Omega(\Gamma_{R \cap T})$.
\end{enumerate}}

\bigskip\noindent
The following lemma is a preliminary to the proof of Theorem 2.2. It is well-known and widely used in the study of Coxeter groups. It can be found, for instance, in the exercises of \cite[Chap. 4]{Bourb1}.

\bigskip\noindent
{\bf Lemma 2.3.} 
{\it
\begin{enumerate}
\item
Let $X,Y$ be two subsets of $S$ and let $w \in W$. Then there exists a unique $(X,Y)$-reduced element lying in the double coset $W_X w W_Y$.
\item
Let $X \subset S$ and $w \in W$. Then $w$ is $(\emptyset,X)$-reduced if and only if $\lg (ws) > \lg (w)$ for all $s \in X$, and $\lg (ws) > \lg (w)$ for all $s \in X$ if and only if $\lg (wu) = \lg(w) + \lg(u)$ for all $u \in W_X$.
\item
Let $X \subset S$ and $w \in W$. Then $w$ is $(X,\emptyset)$-reduced if and only if $\lg (sw) > \lg (w)$ for all $s \in X$, and $\lg (sw) > \lg (w)$ for all $s \in X$ if and only if $\lg (uw) = \lg(u) + \lg(w)$ for all $u \in W_X$.
\end{enumerate}}

\bigskip\noindent
{\bf Proof of Theorem 2.2.}
In order to prove the first part, it suffices to determine a set-map $\pi_T : (W \times \SS^f) \to (W_T \times \SS^f_T)$ which satisfies the following properties.
\begin{itemize}
\item
$\pi_T(u,X)= (u,X)$ for all $(u,X) \in W_T \times \SS^f_T$.
\item
$\pi_T$ is equivariant under (left) action of $W_T$.
\item
If $(u,X) \le (v,Y)$, then $\pi_T(u,X) \le \pi_T(v,Y)$.
\end{itemize}

\bigskip\noindent
Let $(u,X) \in W \times \SS^f$. Write $u=u_0u_1$, where $u_0 \in W_T$ and $u_1$ is $(T,\emptyset)$-reduced. Let $X_0 = T \cap u_1 X u_1^{-1}$. Then we set 
\[
\pi_T(u,X) = (u_0,X_0)\,.
\]
Note that, as $W_{X_0} \subset u_1 W_X u_1^{-1}$, the group $W_{X_0}$ is finite, thus $X_0 \in \SS_T^f$.

\bigskip\noindent
We obviously have $\pi_T(u,X)=(u,X)$ for all $(u,X) \in W_T \times \SS^f_T$, and $\pi_T$ is equivariant under the action of $W_T$. So, it remains to show that, if $(u,X), (v,Y) \in W \times \SS^f$ are such that $(u,X) \le (v,Y)$, then $\pi_T(u,X) \le \pi_T(v,Y)$.

\bigskip\noindent
Let $(u,X), (v,Y) \in W \times \SS^f$ such that $(u,X) \le (v,Y)$. Write $u=u_0u_1$ and $v=v_0v_1$, where $u_0,v_0 \in W_T$ and $u_1,v_1$ are $(T,\emptyset)$-reduced. Set $X_0 = T \cap u_1 W_X u_1^{-1}$ and $Y_0=T \cap v_1 W_Y v_1^{-1}$. Then $\pi_T(u,X) = (u_0,X_0)$ and $\pi_T(v,Y) = (v_0,Y_0)$. Let $w=v^{-1}u$ and $w_0=v_0^{-1}u_0$. Since $(u,X) \le (v,Y)$, we have $X \subset Y$, $w \in W_Y$, and $w$ is $(\emptyset,X)$-reduced. We should show that $X_0 \subset Y_0$, $w_0 \in W_{Y_0}$, and $w_0$ is $(\emptyset,X_0)$-reduced. We argue by induction on the length of $w$. It is easily shown that, if $w=1$, then  $u_0=v_0$ and $X_0 \subset Y_0$, hence $\pi_T(u,X) \le \pi_T(v,Y)$. So, we can assume that $\lg(w) \ge 1$ plus the inductive hypothesis. 

\bigskip\noindent
We write $w=sw'$, where $s \in Y$, $w' \in W_Y$, and $\lg(w') = \lg(w)-1$. Let $v'=vs$. The element $(v')^{-1}u = w'$ lies in $W_Y$ and is $(\emptyset,X)$-reduced (because $w$ is $(\emptyset,X)$-reduced), thus $(u,X) \le (v',Y)$. Write $v'=v_0' v_1'$, where $v_0' \in W_T$ and $v_1'$ is $(T, \emptyset)$-reduced, and set $Y_0' = T \cap (v_1') W_Y (v_1')^{-1}$. By induction hypothesis we have $(u_0,X_0) = \pi_T(u,X) \le \pi_T(v',Y) = (v_0',Y_0')$. Write $w_0' = (v_0')^{-1}u_0$. So, $X_0 \subset Y_0'$, $w_0' \in W_{Y_0'}$, and $w_0'$ is $(\emptyset,X_0)$-reduced. 

\bigskip\noindent
Suppose that $v_1s$ is $(T,\emptyset)$-reduced. Then $v_0'=v_0$ and $v_1'=v_1s$. Moreover, it it easily shown that, in that case, $Y_0 = Y_0'$ (thus $X_0 \subset Y_0$), $w_0=w_0'\in W_{Y_0}$, and $w_0$ is $(\emptyset,X_0)$-reduced. That is, $\pi_T(u,X) \le \pi_T(v,Y)$. 

\bigskip\noindent
So, we can assume that $v_1s$ is not $(T,\emptyset)$-reduced. We have $\lg(v_1s)>\lg(v_1)$, otherwise $v_1s$ would be $(T,\emptyset)$-reduced since $v_1$ is $(T,\emptyset)$-reduced. On the other hand, by Lemma~2.3, there exists $t \in T$ such that $\lg(tv_1s) < \lg (v_1s)$. We also have $\lg(tv_1) > \lg(v_1)$ because $v_1$ is $(T,\emptyset)$-reduced. By the exchange condition (see \cite{Brown1}, Page 47), these inequalities imply that $tv_1=v_1s$. Then we have $v_0'=v_0t$, $v_1'=v_1$, and, therefore, $Y_0=Y_0'$ and $w_0=tw_0'$. A straightforward consequence of this is that $X_0 \subset Y_0=Y_0'$ and $w_0 \in W_{Y_0}$ (since $w_0' \in W_{Y_0}$ and $t = v_1 s (v_1)^{-1} \in T \cap v_1 W_Y v_1^{-1} = Y_0$). It remains to show that $w_0$ is $(\emptyset,X_0)$-reduced. Suppose not. Then we have $\lg(w_0) = \lg(tw_0') > \lg(w_0')$, otherwise $w_0$ would be $(\emptyset,X_0)$-reduced because $w_0'$ is. On the other hand, by Lemma~2.3, there exists $x \in X_0$ such that $\lg(t w_0' x)< \lg(tw_0')$. We also have $\lg( w_0' x) > \lg(w_0')$ since $w_0'$ is $(\emptyset,X_0)$-reduced. By the exchange condition, it follows that $t w_0' = w_0' x = w_0$. Thus:
\[\begin{array}{lc}
&x=(w_0')^{-1} t (w_0') = u_0^{-1} (v_0') t (v_0')^{-1} u_0 \in W_{X_0} = W_T \cap u_1 W_X u_1^{-1}\\
\Rightarrow \quad & u^{-1} v_0 t v_0^{-1} u = u^{-1} v s v^{-1} u = w^{-1}sw \in W_X\\
\Rightarrow \quad & sw W_X = w' W_X = w W_X\,.
\end{array}\]
This contradicts the fact that $w$ is $(\emptyset,X)$-reduced (recall that $\lg(w') < \lg(w)$). So, $w_0$ is $(\emptyset,X_0)$-reduced and, therefore, $\pi_T(u,X) \le \pi_T(v,Y)$.

\bigskip\noindent
In order to prove the second part of the theorem, it suffices to show that the restriction of $\pi_T$ to $(W_R \times \SS^f_R)$ coincides with $\pi_{R \cap T}: (W_R \times \SS^f_R) \to (W_{R \cap T} \times \SS^f_{R \cap T})$. Let $(u,X) \in W_R \times \SS^f_R$. Write $u=u_0u_1$, where $u_0 \in W_T$ and $u_1$ is $(T,\emptyset)$-reduced, and set $X_0 = T \cap u_1 W_X u_1^{-1}$. We have $\pi_T (u,X) = (u_0,X_0)$. Since $u \in W_R$ and $\lg(u) = \lg(u_0) + \lg(u_1)$, we have $u_0 \in W_R$, thus $u_0\in W_R \cap W_T = W_{R \cap T}$. Moreover, by Lemma 2.3, $u_1$ is $(R \cap T, \emptyset)$-reduced (since it is $(T,\emptyset)$-reduced and $R \cap T \subset T$). Finally, $u_1 \in W_R$, since $u \in W_R$ and $\lg(u) = \lg(u_0) + \lg(u_1)$, and $W_X \subset W_R$, since $X \subset R$, thus
\begin{multline*}
X_0 = T \cap u_1 W_X u_1^{-1} = S \cap W_T \cap W_R \cap u_1 W_X u_1^{-1} \\
= S \cap W_{R \cap T} \cap u_1 W_X u_1^{-1} = R \cap T \cap u_1 W_X u_1^{-1}\,.
\end{multline*}
This shows that $\pi_{R \cap T} (u,X) = (u_0,X_0)$.
\qed

\bigskip\noindent
Theorem 2.2 immediately implies the following, which, curiously, was unknown before.

\bigskip\noindent
{\bf Corollary 2.4.}
{\it Let $T$ be a subset of $S$. If $\Gamma$ is of type $K(\pi,1)$, then $\Gamma_T$ is also of type $K(\pi,1)$.}

\bigskip\noindent
Theorem 2.2 also gives a new proof of the following result. (We point out here that Theorem~2.5 is not needed in any known proof of Theorem 2.1!).

\bigskip\noindent
{\bf Theorem 2.5} (Van der Lek \cite{Lek1}).
{\it \begin{enumerate}
\item
Let $T$ be a subset of $S$. Then the embedding $\Sigma_T\hookrightarrow \Sigma$ induces an injective homomorphism $A_{\Gamma_T} \to A_\Gamma$. In other words, $(A_T,\Sigma_T)$ is an Artin-Tits system of $\Gamma_T$.
\item
Let $R$ and $T$ be two subsets of $S$. Then $A_R \cap A_T = A_{R \cap T}$.
\end{enumerate}}

\bigskip\noindent
{\bf Proof.}
Let $T$ be a subset of $S$. We have the following commutative diagram, where the lines are exact sequences.
\[\begin{array}{ccccccccc}
1 & \to  & CA_{\Gamma_T} & \longrightarrow & A_{\Gamma_T} & \longrightarrow & W_{\Gamma_T} & \to & 1\\
  &      & \downarrow    &                 & \downarrow   &                 & \downarrow\\
1 & \to  & CA_{\Gamma}   & \longrightarrow & A_{\Gamma}   & \longrightarrow & W_{\Gamma}   & \to & 1
\end{array}\]
The homomorphism $CA_{\Gamma_T} \to CA_\Gamma$ is injective by Theorem 2.2, and the homomorphism $W_{\Gamma_T} \to W_\Gamma$ is injective by \cite{Bourb1}. We conclude by the five lemma that the homomorphism $A_{\Gamma_T} \to A_\Gamma$ is injective, too.

\bigskip\noindent
For a subset $T \subset S$ we set $CA_T = CA \cap A_T$. Clearly, $CA_T$ is the colored Artin-Tits  group of $\Gamma_T$. Let $R$ and $T$ be two subsets of $S$. In order to prove the second part of the theorem, it suffices to show the inclusion $A_R \cap A_T \subset A_{R \cap T}$: the reverse inclusion $A_{R \cap T} \subset A_R \cap A_T$ is obvious. Let $\alpha \in A_R \cap A_T$. Set $u=\theta(\alpha)$. Note that $u \in W_R \cap W_T= W_{R \cap T}$. Set $\alpha_0=\tau(u)$. By construction, we have $\alpha_0 \in A_{R \cap T}$. Let $\beta = \alpha \alpha_0^{-1}$. We have $\beta \in CA$ and $\beta \in A_R \cap A_T$, thus $\beta \in CA_R \cap CA_T$. We have $\pi_{T\, \ast}(\beta) = \beta$ since $\beta \in CA_T$. On the other hand, by Theorem 2.2, $\pi_{T\, \ast} (\beta) = \pi_{R \cap T\, \ast} (\beta) \in CA_{R \cap T}$ since $\beta \in CA_R$. Thus $\beta \in CA_{R \cap T}$, therefore $\alpha = \beta \alpha_0 \in A_{R \cap T}$.
\qed

\bigskip\noindent
We turn now to combinatorial questions on Artin-Tits groups, and, more precisely, on parabolic subgroups of Artin-Tits groups. Our goal for the remainder of the section is to show that there exists an algorithm which, given a subset $T \subset S$ and a word $\omega \in (\Sigma \sqcup \Sigma^{-1})^\ast$, decides whether the element of $A$ represented by $\omega$ belongs to $A_T$, and, if yes,  determines an expression $\kappa_T(\omega)$ of this element in $(\Sigma_T \sqcup \Sigma_T^{-1})^\ast$ (we will need, of course, to assume that $A$ has a solution to the word problem). We start with an explicit calculation of the homomorphism $\pi_T :CA \to CA_T$ induced by the map $\pi_T: \Omega(\Gamma) \to \Omega (\Gamma_T)$. (Note that, from now on, in order to avoid unwieldy notations, we use the same notation for the map $\Omega(\Gamma) \to \Omega (\Gamma_T)$ as for the induced homomorphism $CA \to CA_T$.)

\bigskip\noindent
For $w \in W$ and $s \in S$ such that $\lg (ws) > \lg(w)$, we set
\[
\delta (w,s) = \tau (w) \sigma_s^2 \tau(w)^{-1}\,,
\]
where $\tau : W \to A$ is the set-section of $\theta : A \to W$ defined in Section 1. It is easily seen that $\delta(w,s) \in CA$ for all $w \in W$ and $s \in S$ such that $\lg (ws) > \lg(w)$, and that the set $\{ \delta(w,s); w \in W,\ s \in S, \text{ and } \lg(ws)>\lg(w) \}$ generates $CA$. Moreover, there exists an algorithm which, given an element $\alpha \in CA$ and an expression $\omega \in (\Sigma \sqcup \Sigma^{-1})^\ast$ of $\alpha$, determines $w_1, \dots ,w_m \in W$, $t_1, \dots, t_m \in S$, and $\mu_1, \dots, \mu_m \in \{ \pm 1\}$, such that $\lg(w_it_i) > \lg(w_i)$ for all $1 \le i \le m$, and
\[
\alpha = \delta(w_1,t_1)^{\mu_1} \cdots \delta(w_m,t_m)^{\mu_m}\,.
\]
The proof of this last fact is left to the reader.

\bigskip\noindent
Let $\Upsilon$ be an abstract simplicial complex. A {\it combinatorial path} in $\Upsilon$ is defined to be a sequence $(v_0,v_1, \dots, v_n)$ of vertices such that $v_{i-1}$ is linked to $v_i$ by an edge for all $1 \le i \le n$. Let $\gamma = (v_0, v_1, \dots, v_n)$ be a combinatorial path in $\Upsilon$. Let $a_i$ be the edge in  the geometric realization $|\Upsilon|$ which links $v_{i-1}$ to $v_i$ and to which we attribute an orientation from $v_{i-1}$ to $v_i$. Then $|\gamma| = a_1 a_2 \cdots a_n$ is a path in $|\Upsilon|$ from $v_0$ to $v_n$. It is called the {\it geometric realization} of $\gamma$.

\bigskip\noindent
We denote by $p: \Omega \to \Omega/W$ the natural projection. We set $x_0 = (1, \emptyset) \in \Omega$ and $\bar x_0 = p(x_0)$. For all $s \in S$ we define the combinatorial path
\[
\gamma_s = ((1,\emptyset), (1, \{s\}), (s,\emptyset))\,.
\]
It is easily deduced from \cite{ChaDav2} and \cite{Salve1} that the natural isomorphism $A \to \pi_1(\Omega/W, \bar x_0)$ sends $\sigma_s$ to $p(|\gamma_s|)$ for all $s \in S$. 

\bigskip\noindent
{\bf Lemma 2.6.}
{\it Let $T$ be a subset of $S$. Let $w \in W$ and $s \in S$ such that $\lg (ws) > \lg (w)$. Write $w=w_0w_1$, where $w_0 \in W_T$ and $w_1$ is $(T,\emptyset)$-reduced. If $w_1s$ is $(T,\emptyset)$-reduced, then $\pi_T(\delta(w,s))=1$. If $w_1s$ is not $(T,\emptyset)$-reduced, then there exists $t \in T$ such that $w_1 s = t w_1$ and $\pi_T(\delta(w,s)) = \delta(w_0,t)$.}

\bigskip\noindent
{\bf Proof.}
For $w \in W$ and $s \in S$ we denote by $\gamma (w,s)$ the combinatorial path
\[
\gamma(w,s) = ((w,\emptyset), (w,\{s\}), (ws, \emptyset))\,.
\]
Note that $\gamma(w,s) = w \cdot \gamma_s$. In order to avoid unwieldy notations, we will still denote by $\gamma(w,s)$ the geometric realization of $\gamma(w,s)$. We start  calculating $\pi_T(\gamma(w,s))$ and $\pi_T(\gamma(ws,s))$, under the assumption that $\lg(ws) > \lg (w)$. Write $w=w_0 w_1$, where $w_0 \in W_T$ and $w_1$ is $(T,\emptyset)$-reduced. We have two cases, depending on whether $w_1 s$ is $(T,\emptyset)$-reduced or not.

\bigskip\noindent
Suppose that $w_1s$ is $(T,\emptyset)$-reduced. Then $T \cap w_1 W_{\{s\}} w_1^{-1} = \emptyset$, that is $w_1 s w_1^{-1} \not\in T$, otherwise it would exist $t \in T$ such that $w_1 s = t w_1$, and this would mean that $w_1 s$ is not $(T,\emptyset)$-reduced. It follows that   
\[
\pi_T((w,\emptyset)) = \pi_T((w,\{s\})) = \pi_T((ws,\emptyset)) = \pi_T((ws,\{s\})) = (w_0,\emptyset)\,.
\] 
So, $\pi_T(\gamma(w,s)) = \pi_T( \gamma(ws,s))$ is the constant path on $(w_0,\emptyset)$.

\bigskip\noindent
Suppose that $w_1s$ is not $(T,\emptyset)$-reduced. Then there exists $t \in T$ such that $\lg(tw_1s) < \lg(w_1s)$. Furthermore, we have $\lg(w_1s) > \lg(w_1)$ and $\lg(tw_1) > \lg(w_1)$ (since $w_1$ is $(T,\emptyset)$-reduced), thus $tw_1 = w_1s$ and $T \cap w_1 W_{\{s\}} w_1^{-1} = T \cap w_1 s W_{\{s\}} s w_1^{-1} = \{t\}$. Hence
\[\begin{array}{c}
\pi_T((w,\emptyset))=(w_0,\emptyset), \ \pi_T((w,\{s\})) = (w_0,\{t\}),\\
\noalign{\smallskip}
\pi_T((ws,\emptyset)) = (w_0t, \emptyset), \ \pi_T((ws, \{s\})) = (w_0t, \{t\})\,.
\end{array}\]
So,
\[
\pi_T(\gamma(w,s))=\gamma(w_0,t)\quad \text{and}\quad \pi_T(\gamma(ws,s)) = \gamma(w_0t,t)\,.
\]

\bigskip\noindent
Now, we assume that $\lg(ws) > \lg (w)$, and we calculate $\pi_T(\delta(w,s))$. Write $w=w_0w_1$, where $w_0 \in W_T$ and $w_1$ is $(T,\emptyset)$-reduced. Let $w_0=s_1 \cdots s_r$ be a reduced form for $w_0$, and let $w_1 = s_{r+1} \cdots s_n$ be a reduced form for $w_1$. For $0 \le i \le n$ we set $u_i = s_1 \cdots s_i$. By the above, we have
\[
\pi_T(\gamma(u_{i-1},s_i)) = \left\{
\begin{array}{ll}
\gamma(u_{i-1},s_i) &\quad \text{if } 1 \le i \le r\,,\\
\cst_{(w_0,\emptyset)} & \quad \text{if } r+1 \le i \le n\,,
\end{array}\right.
\]
where $\cst_{(w_0,\emptyset)}$ denotes the constant path on $(w_0,\emptyset)$. On the other hand,
\[
\delta(w,s) = \gamma(u_0,s_1) \cdots \gamma(u_{n-1},s_n) \gamma(w,s) \gamma(ws,s) \gamma(u_{n-1},s_n)^{-1} \cdots \gamma (u_0,s_1)^{-1}\,.
\]
If $w_1 s$ is $(T, \emptyset)$-reduced, then, by the above,
\[
\pi_T(\delta(w,s)) = \gamma(u_0,s_1) \cdots \gamma(u_{r-1},s_r) \gamma (u_{r-1},s_r)^{-1} \cdots \gamma(u_0,s_1)^{-1} = 1\,.
\]
Suppose that $w_1 s$ is not $(T,\emptyset)$-reduced. As shown before, there exists $t \in T$ such that $w_1 s = t w_1$. Then, by the above,
\begin{multline*}
\pi_T(\delta(w,s)) = \gamma(u_0,s_1) \cdots \gamma(u_{r-1},s_r) \gamma(w_0,t) \gamma(w_0t,t) \gamma (u_{r-1},s_r)^{-1} \cdots \gamma(u_0,s_1)^{-1}
\\
 = \delta(w_0,t)\,.
\qed
\end{multline*}

\bigskip\noindent
There is an algorithm which, given $w \in W$, determines $w_0,w_1 \in W$ such that $w=w_0w_1$, $w_0 \in W_T$, and $w_1$ is $(T,\emptyset)$-reduced. In particular, this algorithm decides whether $w$ is $(T,\emptyset)$-reduced or not. There is also an algorithm which, given $w_1 \in W$ and $s \in S$ such that $w_1$ is $(T,\emptyset)$-reduced and $w_1s$ is not $(T,\emptyset)$-reduced, determines $t \in T$ such that $w_1 s = t w_1$. These algorithms can be easily derived from the classical combinatorial theory of Coxeter groups (see \cite{AbrBro1}, \cite{Bourb1}, \cite{Brown1}, or \cite{Davis1}, for instance). By Lemma 2.6, it follows that there exists an algorithm which, given $\alpha \in CA$, a word $\omega \in (\Sigma \sqcup \Sigma^{-1})^\ast$ which represents $\alpha$, and a subset $T \subset S$, determines a word $\tilde \pi_T(\omega) \in (\Sigma_T \sqcup \Sigma_T^{-1})^\ast$ which represents $\pi_T(\alpha)$. Thanks to this, we can now prove the following.

\bigskip\noindent
{\bf Proposition 2.7.}
{\it Assume that $A$ has a solution to the word problem. Then there exists an algorithm which, given $\alpha \in A$, a word $\omega \in (\Sigma \sqcup \Sigma^{-1})^\ast$ which represents $\alpha$, and $T \subset S$, decides whether $\alpha \in A_T$ and, if yes, determines a word $\kappa_T(\omega) \in (\Sigma_T \sqcup \Sigma_T^{-1})^\ast$ which represents $\alpha$.}

\bigskip\noindent
{\bf Proof.}
Set $w= \theta(\alpha)$. If $w \not \in W_T$, then $\alpha \not \in A_T$. So, we can assume that $w \in W_T$. Recall that there exists an algorithm which determines a reduced form $w=s_1 \cdots s_n$ for $w$ (see \cite{Brown1} for example). Furthermore, we have $w \in W_T$ if and only if $s_1, \dots, s_n \in T$. Set $\tilde \tau(w) = \sigma_{s_1} \cdots \sigma_{s_n}$. Then, by the above, $\tilde \tau(w) \in (\Sigma_T \sqcup \Sigma_T^{-1})^\ast$. Set $\beta = \alpha \tau(w)^{-1}$. We have $\alpha \in A_T$ if and only if $\beta \in CA \cap A_T = CA_T$, we have $\beta \in CA_T$ if and only if $\pi_T(\beta) = \beta$, and we have $\pi_T(\beta) = \beta$ if and only if $\tilde \pi_T(\omega \tilde\tau(w)^{-1})$ and $\omega \tilde\tau(w)^{-1}$ represent the same element in $A$. Applying a solution to the word problem in $A$, we can check whether $\tilde \pi_T(\omega \tilde\tau(w)^{-1})$ and $\omega \tilde\tau(w)^{-1}$ represent the same element in $A$. Moreover, if they represent the same element, then $\kappa_T(\omega) = \tilde \pi_T(\omega \tilde\tau(w)^{-1}) \tilde\tau(w)$ is a word in $(\Sigma_T \sqcup \Sigma_T^{-1})^\ast$ which represents $\alpha$.
\qed




\section{Charney-Davis-Deligne Complex}

The main tool in the study of the $K(\pi,1)$ problem in $\cite{ChaDav1}$ is a simplicial complex (which admits a cubical structure, but this will be seen later) having the same homotopy type as the universal cover of the space $E$ defined in Section 1. This complex is a CAT(0) cube complex if (and only if) the Coxeter graph is of FC type. It is also CAT(0) if the Coxeter graph is so-called of dimension 2, but with a metric different from its cubical metric. By standard arguments, in both cases, the existence of a CAT(0) metric on $\Phi$ implies that $\Phi$ is contractible and, therefore, that $E$ is an Eilenberg MacLane space. This complex is defined as follows.

\bigskip\noindent
Recall that $\SS^f$ denotes the set of $X \subset S$ such that $W_X$ is finite. We denote by $\AA = \AA(\Gamma,\SS^f)$ the set of cosets $\alpha A_X$ with $\alpha \in A$ and $X \in \SS^f$, which we order by inclusion. Then the {\it Charney-Davis-Deligne complex} is $\Phi = \Phi(\Gamma,\SS^f)=|\AA'|$, the geometric realization of the derived complex of $\AA$. 

\bigskip\noindent
In order to prove that $\Phi$ has the same homotopy type as the universal cover of $E$, Charney and Davis strongly use some Deligne's result \cite{Delig1} which says that all spherical type Coxeter graphs are of type $K(\pi,1)$. Our idea in the present section is to construct such a space  $\Phi$ replacing the family $\SS^f$ by some family $\SS$ of subsets of $S$ satisfying: (1) If $X \in \SS$ and $Y \subset X$, then $Y \in \SS$, (2) $\Gamma_X$ is of type $K(\pi,1)$ for all $X \in \SS$. For some technical reason, we should also add the condition: (3) $X \in \SS$ if $\Gamma_X$ is of spherical type (otherwise our space $\Phi$ would not have the same homotopy type as the universal cover of $E$). A family $\SS$ of subsets of $S$ which satisfy (1), (2) and (3) will be called {\it complete and $K(\pi,1)$}. 

\bigskip\noindent
Let $\SS$  be a complete and $K(\pi,1)$ family of subsets of $S$. We denote by $\AA = \AA(\Gamma,\SS)$ the set of cosets $\alpha A_X$ with $\alpha \in A$ and $X \in \SS$, which we order by inclusion. Then the {\it Charney-Davis-Deligne complex} relative to $\SS$ is $\Phi = \Phi(\Gamma,\SS)=|\AA'|$, the geometric realization of the derived complex of $\AA$. The main result of the present section is the following.

\bigskip\noindent
{\bf Theorem 3.1.}
{\it Let $\SS$ be a complete and $K(\pi,1)$ family of subsets of $S$. Then $\Phi= \Phi(\Gamma,\SS)$ has the same homotopy type as the universal cover of $\Omega = \Omega(\Gamma)$.} 

\bigskip\noindent
Our proof of Theorem 3.1 is independent from the proof of Charney and Davis (see \cite{ChaDav1}) for the case $\SS = \SS^f$. Actually, we do not know how to extend the proof of Charney and Davis to the general case. 

\bigskip\noindent
We start with the description of the universal cover of $\Omega$. We define a relation $\le$ on $A \times \SS^f$ by: $(\alpha,X) \le (\beta,Y)$ if $X \subset Y$ and $\beta^{-1}\alpha$ is of the form $\beta^{-1}\alpha = \tau(w)$, where $w\in W_Y$ and $w$ is $(\emptyset,X)$-reduced. It is easily verified that this relation is an order relation. Then we denote by $\tilde \Omega = \tilde \Omega(\Gamma)$ the geometric realization of the derived complex of $(A \times \SS^f, \le)$.

\bigskip\noindent
{\bf Lemma 3.2.}
{\it We have $\tilde \Omega / CA = \Omega$. In particular, since $\pi_1(\Omega)=CA$, $\tilde \Omega$ is the universal cover of $\Omega$.}
 
\bigskip\noindent
{\bf Proof.}
Observe that the quotient of $A \times \SS^f$ by $CA$ is equal to $W \times \SS^f$ as a set. Let $p: (A \times \SS^f) \to (W \times \SS^f)$ be the quotient map. Then, in order to show the equality $\tilde \Omega /CA = \Omega$, it suffices to prove the following three claims.
\begin{itemize}
\item
Let $(\alpha,X), (\beta,Y) \in A \times \SS^f$. If $(\alpha,X) \le (\beta,Y)$, then $p(\alpha,X) \le p(\beta,Y)$.
\item
Let $(\alpha,X) \in A \times \SS^f$ and $(v,Y) \in W \times \SS^f$. If $p(\alpha,X) \le (v,Y)$, then there exists a unique element $(\beta,Y) \in A \times \SS^f$  such that $(\alpha,X) \le (\beta,Y)$ and $p(\beta,Y) = (v,Y)$.
\item
Let $(\beta,Y) \in A \times \SS^f$ and $(u,X) \in W \times \SS^f$. If $(u,X) \le p(\beta,Y)$, then there exists a unique $(\alpha,X) \in A \times \SS^f$ such that $(\alpha,X) \le (\beta,Y)$ and $p(\alpha,X) = (u,X)$.
\end{itemize}
The proofs of these three claims are left to the reader.
\qed

\bigskip\noindent
Now, the proof of Theorem 3.1 is a direct application to $\tilde \Omega$ of the following result. This is a straightforward and well-known consequence of some Weil's classical result (see \cite{Weil1}). 

\bigskip\noindent
{\bf Theorem 3.3.}
{\it Let $\Upsilon$ be a simplicial complex. Let $\UU$ be a family of simplicial subcomplexes of $\Upsilon$ satisfying the following properties.
\begin{itemize}
\item
$\cup_{U \in \UU}U = \Upsilon$.
\item
For all $x \in \Upsilon$, the set of $U \in \UU$ such that $x \in U$ is finite (we say that $\UU$ is {\rm locally finite}).
\item
Let $U_1,U_2 \in \UU$. If $U_1 \cap U_2 \neq \emptyset$, then $U_1 \cap U_2 \in \UU$.
\item
Every $U \in \UU$ is contractible.
\end{itemize}
We order $\UU$ by the inclusion, and we denote by $\Phi= |\UU'|$ the geometric realization of the derived complex of $\UU$. Then $\Phi$ has the same homotopy type as  $\Upsilon$.}

\bigskip\noindent
{\bf Proof of Theorem 3.1.}
Let $\Upsilon$ be a simplicial complex, and let $X$ be a set of vertices of $\Upsilon$. Then the {\it full subcomplex of $\Upsilon$ generated by $X$} is defined to be the simplicial complex over $X$ formed by all the simplices of $\Upsilon$ having their vertices in $X$.

\bigskip\noindent
Let $X \in \SS$. The inclusion $A_X \times \SS^f_X \hookrightarrow A \times \SS^f$ induces an embedding $\tilde \Omega (\Gamma_X) \hookrightarrow \tilde \Omega$. Let $\alpha \in A$ and $X \in \SS$. To the coset $\alpha A_X$ one can associate the subcomplex $U(\alpha A_X) = \alpha \tilde \Omega(\Gamma_X)$ of $\tilde \Omega$. Note that $U(\alpha A_X)$ is the full subcomplex of $\tilde \Omega$ generated by $\alpha A_X \times \SS^f_X$. In particular, the definition of $U(\alpha A_X)$ does not depend on the choice of $\alpha \in \alpha A_X$, and we have $U(\alpha A_X) = U(\beta A_Y)$ if and only if $\alpha A_X = \beta A_Y$. Now, in order to apply Theorem 3.3, we need to show the following:
\begin{enumerate}
\item
Let $\alpha A_X, \beta A_Y\in \AA$. We have $U(\alpha A_X) \subset U(\beta A_Y)$ if and only if $\alpha A_X \subset \beta A_Y$.
\item
$\cup_{\alpha A_X \in \AA} U(\alpha A_X) = \tilde \Omega$.
\item
Let $(\alpha,X) \in A \times \SS^f$ be a vertex of $\tilde \Omega$. Then the set of $\beta A_Y \in \AA$ satisfying $(\alpha,X) \in U(\beta A_Y)$ is finite.
\item
Let $\alpha A_X, \beta A_Y \in \AA$ such that $U(\alpha A_X) \cap U(\beta A_Y) \neq \emptyset$. Then there exists $\gamma A_Z \in \AA$ such that $U(\alpha A_X) \cap U(\beta A_Y) = U(\gamma A_Z)$.
\item
$U(\alpha A_X)$ is contractible for all $\alpha A_X \in \AA$.
\end{enumerate}

\bigskip\noindent
We start proving (1). Observe that, if $X$ is a subset of $S$, then $A_X \cap \Sigma = \Sigma_X$. The inclusion $\Sigma_X \subset A_X \cap \Sigma$ is obvious. Reciprocally, if $\sigma_s \in A_X \cap \Sigma$, then $\theta(\sigma_s) = s  \in W_X \cap S = X$, thus $\sigma_s \in \Sigma_X$. Now, let $\alpha A_X, \beta A_Y \in \AA$ such that $U(\alpha A_X) \subset U(\beta A_Y)$. Since $U(\alpha A_X)$ is the full subcomplex generated by $\alpha A_X \times \SS^f_X$ and $U(\beta A_Y)$ is the full subcomplex generated by $\beta A_Y \times \SS^f_Y$, we have $\alpha A_X \subset \beta A_Y$. Suppose that $\alpha A_X \subset \beta A_Y$. Then $\alpha A_Y = \beta A_Y$, thus
\[
\Sigma_X = A_X \cap \Sigma = \alpha^{-1} \alpha A_X \cap \Sigma \subset \alpha^{-1} \alpha A_Y \cap \Sigma = A_Y \cap \Sigma = \Sigma_Y\,,
\]
hence $X \subset Y$, therefore $\SS^f_X \subset \SS^f_Y$. It follows that $(\alpha A_X \times \SS^f_X) \subset (\beta A_Y \times \SS^f_Y)$, thus $U(\alpha A_X) \subset U(\beta A_Y)$.

\bigskip\noindent
Let $(\alpha_0,X_0) < (\alpha_1,X_1) < \cdots < (\alpha_p,X_p)$ be a chain in $A \times \SS^f$, and let $\Delta = \Delta ((\alpha_0,X_0), \dots, 
\break
(\alpha_p,X_p))$ be its associated simplex in $\tilde \Omega$. Then $X_p \in \SS^f \subset \SS$ and $\Delta \subset \alpha_p \tilde \Omega(\Gamma_{X_p}) = U(\alpha_p A_{X_p})$. This shows that $\cup_{\alpha A_X \in \AA} U(\alpha A_X) = \tilde \Omega$. 

\bigskip\noindent
Let $(\alpha, X) \in A \times \SS^f$. Then the set of $\beta A_Y \in \AA$ satisfying $(\alpha,X) \in U(\beta A_Y)$ is $\{\alpha A_Y; X \subset Y \text{ and } Y \in \SS\}$, which is clearly finite.

\bigskip\noindent
Let $\alpha A_X, \beta A_Y \in \AA$ such that $U(\alpha A_X) \cap U(\beta A_Y) \neq \emptyset$. The set of vertices of $U(\alpha A_X)$ being $\alpha A_X \times \SS^f_X$ and the set of vertices of $U(\beta A_Y)$ being $\beta A_Y \times \SS^f_Y$, we have $\alpha A_X \cap \beta A_Y \neq \emptyset$. Take $\gamma \in \alpha A_X \cap \beta A_Y$ and set $Z=X\cap Y$. We turn now to show that $U(\alpha A_X) \cap U(\beta A_Y) = U(\gamma A_Z)$. Note that $\alpha A_X = \gamma A_X$ and $\beta A_Y =\gamma A_Y$, thus the inclusion $U(\gamma A_Z) \subset U(\alpha A_X) \cap U(\beta A_Y)$ is obvious. So, we just need to prove $U(\alpha A_X) \cap U(\beta A_Y) \subset U(\gamma A_Z)$. Since $U(\alpha A_X)$ and $U(\beta A_Y)$ are full subcomplexes of $\tilde \Omega$, $U(\alpha A_X) \cap U(\beta A_Y)$ is also a full subcomplex of $\tilde \Omega$. So, in order to prove $U(\alpha A_X) \cap U(\beta A_Y) \subset U(\gamma A_Z)$, it suffices to show that every vertex of $U(\alpha A_X) \cap U(\beta A_Y)$ is a vertex of $U(\gamma A_Z)$. Let $(\delta,T)$ be a vertex of $U(\alpha A_X) \cap U(\beta A_Y)$. Then $\delta \in \alpha A_X \cap \beta A_Y= \gamma A_X \cap \gamma A_Y = \gamma (A_X \cap A_Y) = \gamma A_Z$, and $T \in \SS^f_X \cap \SS^f_Y = \SS^f_Z$, thus $(\delta,T)$ is a vertex of $U(\gamma A_Z)$.

\bigskip\noindent
Finally, the hypothesis ``$\SS$ is complete and $K(\pi,1)$'' implies that $U(\alpha A_X) = \alpha \tilde \Omega(\Gamma_X)$ is contractible for every coset $\alpha A_X \in \AA$.
\qed




\section{Cubical structure on $\Phi$ and metric properties}

A {\it geodesic segment} in a metric space $(E,d)$ is an isometric embedding of an interval $[0,l]$ into $E$. A metric space $(E,d)$ is called {\it geodesic} if any two points of $E$ are joined by a geodesic segment. A {\it geodesic triangle} $\TT$ is three geodesic segments $\gamma_1, \gamma_2, \gamma_3$ that join three  points in $E$. Then a {\it comparison triangle} for $\TT$ in the Euclidean plane $\E^2$ is defined to be a geodesic triangle $\bar \TT = (\bar \gamma_1, \bar \gamma_2, \bar \gamma_3)$ in the plane $\E^2$ such that the lengths of $\gamma_i$ and $\bar\gamma_i$ are equal for all $i=1,2,3$. Note that such a comparison triangle always exists and is unique up to isometry. For a point $x$ in $\gamma_i$ we denote by $\bar x$ the corresponding point in $\bar \gamma_i$. We say that the triangle $\TT$ is CAT(0) if, for every pair $x,y$ of points in $\TT$, the distance between $x$ and $y$ is less or equal to the distance between $\bar x$ and $\bar y$. A geodesic metric space $(E,d)$ is called CAT(0) if all its triangles are CAT(0). The study of CAT(0) spaces and the discrete groups acting on them is an active domain in mathematics. We refer to \cite{BriHae1} for a complete and detailed account on this theory. We just recall here that CAT(0) spaces are contractible. 

\bigskip\noindent
A {\it cube complex} is a polyhedral cell complex $E$ such that each cell is isometric to some standard cube $[0,1]^n$ in a Euclidean space, and such that the gluing maps are isometries. We should specify here that $E$ must be regular, in the sense that two different facets of a cube cannot be glued together. If the dimension of the cubes is bounded, then such a space is a complete geodesic metric space (see \cite{Brids1}).

\bigskip\noindent
Let $E$ be a cube complex, and let $v$ be a vertex of $E$. The {\it link} of $E$ in $v$ is the set $\link(v,E) = \{x \in E; d(v,x)=\varepsilon\}$ of points in $E$ at distance $\varepsilon$ from $v$, where $\varepsilon$ is a real number strictly less than $\frac{1}{2}$. The cellular decomposition of $E$ induces a cellular decomposition of $\link(v,E)$, where each cell is naturally a simplex, but this decomposition is not necessarily a triangulation (in the sense that this is not necessarily a simplicial complex). For example, if we glue two squares along their boundaries, then the link of the space at any vertex is not a simplicial complex. We say that $E$ is {\it locally regular} if $\link(v,E)$ is a simplicial complex for every vertex $v$ of $E$.   

\bigskip\noindent
Let $E$ be a locally regular cube complex, and let $v$ be a vertex of $E$. Then $\link(v,E)$ is the geometric realization of an abstract simplicial complex  $\Link(v,E)$ defined as follows. (1) The vertices of $\Link(v,E)$ are the vertices $u$ of $E$ such that $[u,v]$ is an edge of $E$; (2) $\{u_0, \dots, u_p\}$ is a simplex of $\Link(v,E)$ if there exists a cube (cell) in $E$ containing $\{v,u_0,u_1, \dots, u_p\}$ in its set of vertices.    

\bigskip\noindent
Let $\Upsilon$ be an (abstract) simplicial complex, and let $V$ be its set of vertices. A subset $\Delta$ of $V$ is called a {\it presimplex} if $\{u,v\}$  is an edge ({\it i.e.} a 1-simplex) for all $u,v \in \Delta$, $u \neq v$. We say that $\Upsilon$ is a {\it flag complex} if every presimplex is a simplex.

\bigskip\noindent
The following is a simple a useful criterion for a cube complex to be CAT(0).

\bigskip\noindent
{\bf Theorem 4.1} (Gromov \cite{Gromo1}).
{\it Let $E$ be a simply connected finite dimensional cube complex. Then $E$ is CAT(0) if and only if it is locally regular and $\Link(v,E)$ is a flag complex for every vertex $v$ of $E$.}

\bigskip\noindent
Now, we turn back to the situation where $\Gamma$ is a Coxeter graph, $(W,S)$ is its associated Coxeter system, and $(A, \Sigma) = (A_\Gamma,\Sigma)$ is its associated Artin-Tits system. We assume given a complete and $K(\pi,1)$ family $\SS$ of subsets of $S$, and we set $\Phi = \Phi(\Gamma,\SS)$. Note that this family, being closed under inclusion, can (and will) be considered as a simplicial complex over $S$.     

\bigskip\noindent
We start endowing $\Phi$ with a structure of locally regular cube complex.

\bigskip\noindent 
We set $S = \{s_1, \dots, s_n\}$. Let $K^n=[0,1]^n$ denote the standard cube of dimension $n$, and let $\PP(S)$ denote the set of subsets of $S$. There is a bijection $v$ from $\PP(S)$ to the set of vertices of $K^n$ defined as follows. Let $T \in \PP(S)$. Then $v(T) = (\varepsilon_1, \dots, \varepsilon_n)$, where $\varepsilon_i = 1$ if $s_i \in T$, and $\varepsilon_i = 0$ if $s_i \not\in T$. With every pair $(R,T)$ of subsets of $S$ such that $R \subset T$, one can associate a face of $K^n$, denoted by $K(R,T)$, whose vertices are of the form $v(U)$ with $R \subset U \subset T$. It is easily seen that every face of $K^n$ is of this form. Assume $\PP(S)$ to be endowed with the inclusion relation. Then the geometric realization $|\PP(S)'|$ of the derived complex of $\PP(S)$ is a simplicial decomposition of $K^n$. For $R \subset T \subset S$, the face $K(R,T)$ is the union of the simplices $\Delta(U_0,U_1, \dots, U_p)$ of $|\PP(S)'|$ such that $U_0=R \subset U_1 \subset \cdots \subset U_p=T$.

\bigskip\noindent
For $\alpha \in A$ and $R,T \in \SS$ such that $R \subset T$, there is a map from $K(R,T)$ to $\Phi$ which sends a simplex $\Delta(R,U_1, \dots, U_{p-1}, T)$ to the simplex $\Delta(\alpha A_R, \alpha A_{U_1}, \dots, \alpha A_{U_{p-1}}, \alpha A_T)$ of $\Phi$. It is easily checked that this map is injective and that its image is a cube which depends only on the cosets $\alpha A_R$ and $\alpha A_T$. This cube will be denote by $C(\alpha A_R, \alpha A_T)$. Note that every face of $C(\alpha A_R, \alpha A_T)$ is of the form $C(\alpha A_{U_1}, \alpha A_{U_2})$, where $R \subset U_1 \subset U_2 \subset T$. We leave to the reader to check that the set of cubes of the form $C(\alpha A_R, \alpha A_T)$, with $\alpha \in A$ and $R,S \in \SS$, $R \subset T$, endows $\Phi$ with a structure of locally regular cube complex. From now on, we will always assume $\Phi$ to be endowed with this cubical structure.

\bigskip\noindent
The main result of this section is:

\bigskip\noindent
{\bf Theorem 4.2.}
{\it $\Phi$ is CAT(0) if and only if $\SS$ is a flag complex.}

\bigskip\noindent
{\bf Remark.}
Recall that $\SS_{<\infty}$ denotes the set of subsets $X \subset S$ such that $\Gamma_X$ is free of infinity. Let $\SS$ be a set of subsets of $S$ which is complete and $K(\pi,1)$, and which is a flag complex. Let $X \in \SS_{<\infty}$. We have $\{s,t\} \in \SS^f \subset \SS$ for all $s,t \in X$, $s \neq t$, thus $X\in \SS$. So, $\SS_{<\infty} \subset \SS$. It is easily seen that $\SS_{<\infty}$ is closed under inclusion, contains $\SS^f$, and is a flag complex. Moreover, as $\SS_{<\infty} \subset \SS$, $\Gamma_X$ is of type $K(\pi,1)$ for all $X \in \SS_{<\infty}$. So, the construction of the complex $\Phi = \Phi(\Gamma,\SS)$ is interesting in its own for any (complete and $K(\pi,1)$) family of subsets of $S$, but Theorem 4.2 is useful essentially in the case $\SS = \SS_{<\infty}$. However, it gives a new proof to the following result due to Ellis and Sköldberg.

\bigskip\noindent
{\bf Corollary 4.3} (Ellis, Sköldberg \cite{EllSko1}).
{\it $\Gamma$ is of type $K(\pi,1)$ if and only if $\Gamma_X$ is of type $K(\pi,1)$ for all $X \in \SS_{<\infty}$.}

\bigskip\noindent
The case $\SS = \SS^f$ in Theorem 4.2 is one of the main results in \cite{ChaDav1}. Actually, our proof is substantially the same as the proof in \cite{ChaDav1} except in one point. As in \cite{ChaDav1}, we will show that, for every vertex $v$ of $\Phi$, $\Link(v,\Phi)$ can be decomposed as the join of a finite simplicial complex with some Artin-Tits complex (see Proposition 4.9). This ``Artin-Tits complex'' is always a flag complex (see Proposition 4.5), thus $\Link(v,\Phi)$ is a flag complex if and only if the finite part in the decomposition is a flag complex. In their work \cite{ChaDav1}, Charney and Davis need to consider only Artin-Tits complexes of spherical type, and, in order to prove that these complexes are flag complexes, strongly use the combinatorial study of spherical type Artin-Tits groups made by the first author in \cite{Charn1} and \cite{Charn2}. This study cannot be extended to the other Artin-Tits groups, and our proof of Proposition 4.5 use different tools such as those introduced in Section 2.

\bigskip\noindent
We start the proof of Theorem 4.2 with the study of what we call ``Artin-Tits complexes''.

\bigskip\noindent
For $s \in S$, we set $W^s=W_{S\setminus \{s\}}$. The {\it Coxeter complex} of $\Gamma$, denoted by $\Cox = \Cox(\Gamma)$, is the (abstract) simplicial complex defined by the following data. (1) The set of vertices of $\Cox$ is the set of cosets $\{wW^s; s \in S, w \in W\}$. (2) A family  $\{w_0W^{s_0}, w_1W^{s_1}, \dots, w_pW^{s_p}\}$ is a simplex of $\Cox$ if the intersection $w_0W^{s_0} \cap w_1W^{s_1} \cap \cdots \cap w_pW^{s_p}$ is nonempty. For $X \subset S$, $X \neq \emptyset$, we set $W^X = W_{S \setminus X}$. Let $X \subset S$, $X \neq \emptyset$, and $w \in W$. With the coset $wW^X$ we can associate the simplex $\Delta(wW^X)= \{wW^s; s \in X\}$ of $\Cox$. Every simplex of $\Cox$ is of this form, and we have  $\Delta(uW^X) = \Delta(vW^Y)$ if and only if $uW^X = vW^Y$. We refer to \cite{AbrBro1} for a detailed study on Coxeter complexes. We just mention here that these complexes are flag complexes (see \cite{AbrBro1}, Exercise 3.116). This will be used in the proof of Proposition 4.5.

\bigskip\noindent
For $s \in S$, we set $A^s=A_{S\setminus \{s\}}$. The {\it Artin-Tits complex} of $\Gamma$, denoted by $\Art = \Art(\Gamma)$, is the (abstract) simplicial complex defined by the following data. (1) The set of vertices of $\Art$ is the set of cosets $\{\alpha A^s; s \in S, \alpha \in A\}$. (2) A family $\{\alpha_0 A^{s_0}, \alpha_1 A^{s_1}, \dots, \alpha_p A^{s_p}\}$ is a simplex of $\Art$ if the intersection $\alpha_0 A^{s_0} \cap \alpha_1 A^{s_1} \cap \cdots \cap \alpha_p A^{s_p}$ is nonempty.

\bigskip\noindent
For $X \subset S$, $X \neq \emptyset$, we set $A^X = A_{S \setminus X}$. Let $X \subset S$, $X \neq \emptyset$, and $\alpha \in A$. With the coset $\alpha A^X$ we can associate the simplex $\Delta(\alpha A^X)= \{\alpha A^s; s \in X\}$ of $\Art$. As for the Coxeter complexes, we have the following.

\bigskip\noindent
{\bf Lemma 4.4.}
{\it If $\Delta$ is a simplex of $\Art$, then there exist $\alpha \in A$ and $X \subset S$, $X \neq \emptyset$, such that $\Delta = \Delta(\alpha A^X)$. Moreover, we have $\Delta(\alpha A^X) = \Delta(\beta A^Y)$ if and only if $\alpha A^X = \beta A^Y$.}

\bigskip\noindent
{\bf Proof.} 
Let $\Delta$ be a simplex of $\Art$. Let $\alpha_0 A^{s_0}, \dots, \alpha_p A^{s_p}$ be the vertices of $\Delta$. By definition, the intersection of the $\alpha_i A^{s_i}$'s is nonempty. Choose some $\beta$ in $\cap_{i=0}^p \alpha_i A^{s_i}$. Then the vertices of $\Delta$ are $\beta A^{s_0}, \dots, \beta A^{s_p}$. Let $i,j \in \{0, \dots, p\}$, $i\neq j$. Since $\beta A^{s_i}$ is different from $\beta A^{s_j}$, we have $s_i \neq s_j$. Set $X=\{s_0, \dots, s_p\}$. Then $X$ is of cardinality $p+1$ and $\Delta = \Delta(\beta A^X)$.

\bigskip\noindent
Let  $\alpha \in A$ and $X\subset S$, $X \neq \emptyset$. We have $\cap_{s \in X} A^s = A^X$ (see Theorem 2.5), thus $\cap_{s \in X} \alpha A^s = \alpha A^X$. Since $\{\alpha A^s ; s\in X\}$ is the set of vertices of $\Delta(\alpha A^X)$, this shows that $\Delta(\alpha A^X)$ determines $\alpha A^X$.
\qed

\bigskip\noindent
{\bf Proposition 4.5.}
{\it $\Art$ is a flag complex.}

\bigskip\noindent
The following lemmas 4.6 and 4.7 are preliminaries to the proof of Proposition 4.5.

\bigskip\noindent
{\bf Lemma 4.6.}
{\it Let $X,Y$ be two subsets of $S$, and let $\alpha \in CA$. If $CA_X \cap \alpha CA_Y \neq \emptyset$, then $\pi_X (\alpha) \in CA_X \cap \alpha CA_Y$.}

\bigskip\noindent
{\bf Proof.}
Let $\beta \in CA_X \cap \alpha CA_Y$. Write $\alpha= \beta \gamma$, where $\gamma \in CA_Y$. By Theorem 2.2, we have $\pi_X (\gamma) = \pi_{X \cap Y}(\gamma)$. In particular, $\pi_X(\gamma) \in CA_{X \cap Y} \subset CA_Y$. Thus,
\[
\pi_X(\alpha) = \pi_X (\beta) \pi_X (\gamma) = \beta \pi_X (\gamma) \in \beta CA_Y = \alpha CA_Y\,.
\qed
\]

\bigskip\noindent
{\bf Lemma 4.7.}
{\it Let $\alpha_1,\alpha_2,\alpha_3 \in A$, and let $X_1,X_2,X_3$ be three subsets of $S$. If $\alpha_i A_{X_i} \cap \alpha_j A_{X_j}\neq\emptyset$ for all $i,j \in \{1,2,3\}$, then $\alpha_1 A_{X_1} \cap \alpha_2 A_{X_2} \cap \alpha_3 A_{X_3} \neq \emptyset$.}

\bigskip\noindent
{\bf Proof.}
Let $\beta_1, \beta_2, \beta_3 \in CA$, and let $X_1,X_2,X_3$ be three subsets of $S$ such that $\beta_i CA_{X_i} \cap \beta_j CA_{X_j} \neq \emptyset$ for all $i,j \in \{1,2,3\}$. Let $\gamma \in \beta_1^{-1} \beta_2 CA_{X_2} \cap \beta_1^{-1} \beta_3 CA_{X_3}$. Since $CA_{X_1} \cap \beta_1^{-1} \beta_2 CA_{X_2} \neq \emptyset$ and $CA_{X_1} \cap \beta_1^{-1} \beta_3 CA_{X_3} \neq \emptyset$, by Lemma 4.6, we have $\pi_{X_1} (\gamma) \in CA_{X_1} \cap \beta_1^{-1} \beta_2 CA_{X_2} \cap \beta_1^{-1} \beta_3 CA_{X_3}$, thus $\beta_1 CA_{X_1} \cap \beta_2 CA_{X_2} \cap \beta_3 CA_{X_3} \neq \emptyset$.

\bigskip\noindent
Now, let $\alpha_1,\alpha_2,\alpha_3 \in A$, and let $X_1,X_2,X_3$ be three subsets of $S$ such that $\alpha_i A_{X_i} \cap \alpha_j A_{X_j}\neq\emptyset$ for all $i,j \in \{1,2,3\}$. Set $u_i = \theta(\alpha_i)$ for all $i=1,2,3$. We have $u_i W_{X_i} \cap u_j W_{X_j} \neq \emptyset$ for all $i,j \in \{ 1,2,3\}$, thus, since $\Cox$ is a flag complex, the intersection $u_1 W_{X_1} \cap u_2 W_{X_2} \cap u_3 W_{X_3}$ is nonempty. Choose $w \in u_1 W_{X_1} \cap u_2 W_{X_2} \cap u_3 W_{X_3}$. Without loss of generality, we can assume that $\theta(\alpha_i) = w$ ({\it i.e.} $u_i=w$) for all $i =1,2,3$. Set $\beta_i = \tau(w)^{-1} \alpha_i$ for all $i=1,2,3$. Note that $\beta_i \in CA$ for all $i \in \{1,2,3\}$, and $\beta_i A_{X_i} \cap \beta_j A_{X_j}\neq\emptyset$ for all $i,j \in \{1,2,3\}$. Now, we show that $\beta_i CA_{X_i} \cap \beta_j CA_{X_j} \neq \emptyset$ for all $i,j \in \{1,2,3\}$. By the previous observation, this implies that $\cap_{i=1}^3 \beta_i CA_{X_i} \neq \emptyset$, thus $\cap_{i=1}^3 \beta_i A_{X_i} \neq \emptyset$, therefore $\cap_{i=1}^3 \alpha_i A_{X_i} \neq \emptyset$.

\bigskip\noindent
Let $i,j \in \{1,2,3\}$, $i \neq j$. Let $\gamma \in \beta_i A_{X_i} \cap \beta_j A_{X_j}$. We write $\gamma = \beta_i \beta_i' = \beta_j \beta_j'$, where $\beta_i' \in A_{X_i}$ and $\beta_j' \in A_{X_j}$. Since $\beta_i, \beta_j \in CA$, we have $\theta(\beta_i') = \theta(\beta_j') \in W_{X_i} \cap W_{X_j} = W_{X_i \cap X_j}$. Let $w'$ be this element. Then $\tau(w') \in A_{X_i \cap X_j} = A_{X_i} \cap A_{X_j}$, thus  
\[
\beta_i \beta_i' \tau(w')^{-1} = \beta_j \beta_j' \tau(w')^{-1} \in \beta_i CA_{X_i} \cap \beta_j CA_{X_j}\,.
\qed
\]

\bigskip\noindent 
{\bf Proof of Proposition 4.5.}
Let $U= \{\alpha_0 A^{s_0}, \alpha_1 A^{s_1}, \dots, \alpha_p A^{s_p} \}$ be a family of vertices of $\Art$ such that $\alpha_i A^{s_i} \cap \alpha_j A^{s_j} \neq \emptyset$ for all $i,j \in \{ 0,1, \dots, p\}$. We show that $\cap_{i=0}^p \alpha_i A^{s_i} \neq \emptyset$ by induction on $p$. The case $p=1$ is trivial and the case $p=2$ is proved in Lemma 4.7, thus we can assume that $p\ge 3$, plus the inductive hypothesis.

\bigskip\noindent
Set $X=\{s_0,s_1, \dots, s_{p-2}\}$. By induction, $\cap_{i=0}^{p-2} \alpha_i A^{s_i} \neq \emptyset$. Let $\beta \in \cap_{i=0}^{p-2} \alpha_i A^{s_i}$. Then $\cap_{i=0}^{p-2} \alpha_i A^{s_i} 
\break
= \beta A^X$. By induction again, $\beta A^X \cap \alpha_{p-1} A^{s_{p-1}} \neq \emptyset$ and $\beta A^X \cap \alpha_p A^{s_p} \neq \emptyset$. By the starting hypothesis, $\alpha_{p-1} A^{s_{p-1}} \cap \alpha_p A^{s_p} \neq \emptyset$. By Lemma 4.7, we conclude that 
\[
\cap_{i=0}^{p} \alpha_i A^{s_i} = \beta A^X \cap \alpha_{p-1} A^{s_{p-1}} \cap \alpha_p A^{s_p} \neq \emptyset\,.
\qed
\]

\bigskip\noindent
Now, we turn back to the study of the Charney-Davis-Deligne complex. 

\bigskip\noindent
Let $\Upsilon_1$ and $\Upsilon_2$ be two simplicial complexes, and let $V_1$ and $V_2$ be their respective sets of vertices. The {\it join}  of $\Upsilon_1$ and $\Upsilon_2$ is the simplicial complex $\Upsilon_1 \ast \Upsilon_2$ defined by: (1) $V_1 \sqcup V_2$ is the set of vertices of $\Upsilon_1 \ast \Upsilon_2$, (2) $\Delta \subset V_1 \sqcup V_2$ is a simplex of $\Upsilon_1 \ast \Upsilon_2$ if $\Delta \cap V_1$ is a simplex of $\Upsilon_1$ and $\Delta \cap V_2$ is a simplex of $\Upsilon_2$. It is clear that $\Upsilon_1 \ast \Upsilon_2$ is a flag complex if and only if $\Upsilon_1$ and $\Upsilon_2$ are both flag complexes.

\bigskip\noindent
Recall that, for $\alpha \in A$ and $R,T \in \SS$, $R \subset T$, we denote by $C(\alpha A_R,\alpha A_T)$ the cube of $\Phi$ union of the simplices $\Delta(\alpha A_{U_0}, \alpha A_{U_1}, \dots, \alpha A_{U_p})$ such that $R=U_0 \subset U_1 \subset \cdots \subset U_p=T$. For $\alpha \in A$ and $T \in \SS$, we set $x(\alpha A_T)=C(\alpha A_T, \alpha A_T)$. This is a vertex of $\Phi$ (viewed as a cube complex), and every vertex of $\Phi$ is of this form. On the other hand, for $T \in \SS$, we denote by $\LL(T) = \LL(T,\SS)$ the simplicial complex defined as follows. The vertices of $\LL(T)$ are the elements $s \in S\setminus T$ such that $T \cup \{s\} \in \SS$. A subset $X \subset S \setminus T$ is a simplex of $\LL(T)$ if $T \cup X \in \SS$. The key point in the proof of Theorem~4.2 is the following. 

\bigskip\noindent
{\bf Proposition 4.9.}
{\it Let $\alpha \in A$ and $T \in \SS$. Then
\[
\Link(x(\alpha A_T),\Phi) \simeq \LL(T) \ast \Art(\Gamma_T)\,.
\]}

\bigskip\noindent
{\bf Proof.}
We denote by $V(\Upsilon)$ the set of vertices of a given simplical complex $\Upsilon$. Consider the map $f: V (\LL(T)) \sqcup V (\Art (\Gamma_T)) \to V (\Link(x(\alpha A_T),\Phi))$ defined as follows. If $s \in V(\LL(T))$, then 
\[
f(s) = x(\alpha A_{T \cup \{s\}})\,.
\]
If $\beta A_T^s= \beta A_{T \setminus \{s\}}$ is a vertex of $\Art(\Gamma_T)$, then
\[
f(\beta A_T^s) = \alpha \beta A_T^s\ = \alpha \beta A_{T \setminus \{s\}}\,.
\]
It is easily seen that $f$ is well-defined and is a one-to-one correspondence. So, it remains to show that $f$ induces a bijection between the set simplices of $\LL(T) \ast \Art(\Gamma_T)$ and the set of simplices of $\Link(x(\alpha A_T),\Phi)$.

\bigskip\noindent
Let $\Delta$ be a simplex of $\LL(T) \ast \Art(\Gamma_T)$. By definition, there exist a simplex $X$ of $\LL(T)$ and a simplex $\Delta (\beta A_T^Y)$ of $\Art(\Gamma_T)$ such that $\Delta = X \sqcup \Delta (\beta A_T^Y)$. Consider the cube $C=C(\alpha \beta A_{T \setminus Y}, 
\break
\alpha A_{T \cup X})$ of $\Phi$. Then $x(\alpha A_T)$ is a vertex of $C$, $f(s)$ is a vertex of $C$ for all $s \in X$, and $f(\beta A_T^t)$ is a vertex of $C$ for all $t \in Y$, thus $f(\Delta)$ is a simplex of $\Link(x(\alpha A_T), \Phi)$.

\bigskip\noindent
Let $\Delta$ be a simplex of $\Link(x(\alpha A_T), \Phi)$. Since $f$ is a bijection, there exist $s_1, \dots, s_p \in V(\LL(T))$ and $\beta_1 A^{t_1}, \dots, \beta_q A^{t_q} \in V(\Art(\Gamma_T))$ such that $f(s_1), \dots, f(s_p), f(\beta_1 A^{t_1}), \dots, f(\beta_q A^{t_q})$ are the vertices of $\Delta$. By definition, there exists a cube $C$ of $\Phi$ such that $x(\alpha A_T)$ is a vertex of $C$, and $f(s_1), \dots, f(s_p), f(\beta_1 A^{t_1}), \dots, f(\beta_q A^{t_q})$ are vertices of $C$. There exist $\gamma \in A$ and $X,Y \in \SS$, $X \subset Y$, such that $C=C(\gamma A_X, \gamma A_Y)$. Set $\Delta_1'= \{s_1, \dots, s_p\}$. For all $1 \le i \le p$ we have
\[
f(s_i) = \alpha A_{T \cup \{s_i\}} \subset \gamma A_Y = \alpha A_Y\,,
\]
thus $s_i \in Y$. This implies that $T \cup \Delta_1' \subset Y$, thus $T \cup \Delta_1' \in \SS$ (since $Y \in \SS$), therefore $\Delta_1'$ is a simplex of $\LL(T)$. Set $\Delta_2' = \{f(\beta_1 A^{t_1}), \dots, f(\beta_q A^{t_q}) \}$. For all $1 \le j \le q$ we have
\[
\alpha^{-1} \gamma A_X \subset \alpha^{-1} f(\beta_j A^{t_j}) = \beta_j A^{t_j}\,,
\]
thus $\cap_{j=1}^q \beta_j A^{t_j} \neq \emptyset$, therefore $\Delta_2'$ is a simplex of $\Art(\Gamma_T)$. Finally, $\Delta' = \Delta_1' \sqcup \Delta_2'$ is a simplex of $\LL(T) \ast \Art(\Gamma_T)$, and $\Delta = f(\Delta')$.
\qed

\bigskip\noindent
{\bf Proof of Theorem 4.2.}
By Theorem 3.1, $\Phi$ is simply connected. So, by Theorem 4.1, $\Phi$ is CAT(0) if and only if $\Link(x(\alpha A_T),\Phi)$ is a flag complex for every vertex $x(\alpha A_T)$ of $\Phi$. By Propositions 4.5 and 4.9, $\Link(x(\alpha A_T),\Phi)$ is a flag complex for every vertex $x(\alpha A_T)$ of $\Phi$ if and only if $\LL(T)$ is a flag complex for every $T \in \SS$. It is easily checked that, for $T \in \SS$, $\LL(T)$ is a flag complex if $\LL(\emptyset)=\SS$ is a flag complex, thus $\LL(T)$ is a flag complex for every $T \in \SS$ if and only if $\SS$ is a flag complex.
\qed




\section{Word problem}

In \cite{AltCha1} Altobelli and Charney use the geometry of $\Phi(\Gamma,\SS^f)$ to determine a solution to the word problem in the FC type Artin-Tits groups. In a similar way, we turn now in this section to use the geometry of $\Phi=\Phi(\Gamma,\SS_{<\infty})$ to determine a solution to the word problem in $A$, when $\Gamma_X$ is of type $K(\pi,1)$ and $A_X$ has a solution to the word problem for all $X \in \SS_{<\infty}$. Our solution is inspired by that of Altobelli and Charney, but is different in the details. Indeed, Altobelli and Charney use several results on spherical type Artin-Tits groups from \cite{Charn1}, \cite{Charn2}, and \cite{Altob1}, and these results cannot be extended to the other Artin-Tits groups. To overcome these difficulties, we will use other tolls including those developed in Section 2.

\bigskip\noindent
Let $E$ be a CAT(0) cube complex. The {\it star} of a cube $C$, denoted by $\Star (C)$, is defined to be the union of all the cubes of $E$ that contain $C$. On the other hand, we say that two cubes $C_1$ and $C_2$ {\it span} a cube if there exists a cube in $E$ containing both, $C_1$ and $C_2$. In that case, the smallest cube of $E$ containing $C_1$ and $C_2$ is called the cube {\it spanned} by $C_1$ and $C_2$, and is denoted by $\Span(C_1,C_2)$.

\bigskip\noindent
Let $x$ and $y$ be two vertices of $E$. A {\it cube path of length $n$} from $x$ to $y$  is a sequence $(C_1, \dots, C_n)$ of $n$ cubes of $E$ such that (1) $x$ is a vertex of $C_1$, and $y$ is a vertex of $C_n$, (2) $C_i \cap C_{i+1} \neq \emptyset$ for all $1 \le i\le n-1$. If $x=y$, we admit there is a (unique) cube path of length 0 from $x$ to $x$. We say that a cube path $\CC = (C_1, \dots, C_n)$ is {\it normal} if (1) $\dim\, C_i \ge 1$ for all $1 \le  i\le n$, (2) $C_i \cap C_{i+1}$ is a vertex, denoted by $x_i$, for all $1 \le i \le n-1$, (3) $C_i=\Span(x_{i-1},x_i)$ for all $1 \le i\le n$, where $x_0=x$ and $x_n=y$, (4) $\Star(C_i) \cap C_{i+1} = \{x_i\}$ for all $1 \le i \le n-1$.

\bigskip\noindent
{\bf Theorem 5.1} (Niblo, Reeves \cite{NibRee1}).
{\it Let $E$ be a CAT(0) cube complex, and let $x$ and $y$ be two vertices of $E$. Then there exists a unique normal cube path from $x$ to $y$ in $E$.}

\bigskip\noindent
Now, let $\Gamma$ be a Coxeter graph, let $(W,S)$ be its associated Coxeter system, and let $(A,\Sigma)$ be its associated Artin-Tits system. We assume that $\Gamma_X$ is of type $K(\pi,1)$ and $A_X$ has a solution to the word problem for all $X \in \SS_{<\infty}$. We set $\Phi = \Phi(\Gamma, \SS_{<\infty})$, which is supposed to be endowed with its cubical structure. So, by Theorem 4.2, $\Phi$ is CAT(0). We start with some technical preliminaries (Lemmas 5.2 -- 5.4).

\bigskip\noindent
{\bf Lemma 5.2.}
{\it Let $C_1 = C(\alpha_1 A_{R_1}, \alpha_1 A_{T_1})$ and $C_2 = C(\alpha_2 A_{R_2}, \alpha_2 A_{T_2})$ be two cubes of $\Phi$. Then:
\begin{enumerate}
\item
$C_1 \cap C_2 \neq \emptyset$ if and only if $R_1 \cup R_2 \subset T_1 \cap T_2$ and $\alpha_1^{-1} \alpha_2 \in A_{T_1 \cap T_2}$. In that case,
\[
C_1 \cap C_2 = C(\alpha_1 A_R, \alpha_1 A_{T_1\cap T_2})\,,
\]
where $R$ is the smallest subset of $S$ satisfying $R_1 \cup R_2 \subset R \subset T_1 \cap T_2$ and $\alpha_1^{-1} \alpha_2 \in A_R$.
\item
$C_1$ and $C_2$ span a cube if and only if $\alpha_1 A_{R_1} \cap \alpha_2 A_{R_2} \neq \emptyset$ and $T_1 \cup T_2 \in \SS_{<\infty}$. In that case,
\[
\Span(C_1,C_2) = C(\alpha A_{R_1 \cap R_2}, \alpha A_{T_1 \cup T_2})\,,
\]
where $\alpha$ is any element in $\alpha_1 A_{R_1} \cap \alpha_2 A_{R_2}$.
\end{enumerate}}

\bigskip\noindent
{\bf Proof.}
Suppose that $C_1 \cap C_2 \neq \emptyset$. Let $x(\beta A_X)$ be a vertex in $C_1 \cap C_2$. Then $\alpha_1 A_{R_1} \subset \beta A_X \subset \alpha_1 A_{T_1}$ and $\alpha_2 A_{R_2} \subset \beta A_X \subset \alpha_2 A_{T_2}$. It follows that $R_1 \subset X \subset T_1$ and $R_2 \subset X \subset T_2$, thus $R_1 \cup R_2 \subset X \subset T_1 \cap T_2$. Moreover, $\alpha_1 A_X = \beta A_X = \alpha_2 A_X$, thus $\alpha_1^{-1} \alpha_2 \in A_X \subset A_{T_1\cap T_2}$.

\bigskip\noindent
Suppose that $R_1 \cup R_2 \subset T_1 \cap T_2$ and $\alpha_1^{-1} \alpha_2 \in A_{T_1\cap T_2}$. Then $\alpha_1 A_{R_1} \subset \alpha_1 A_{T_1 \cap T_2} \subset \alpha_1 A_{T_1}$ and $\alpha_2 A_{R_2} \subset \alpha_2 A_{T_1 \cap T_2} = \alpha_1 A_{T_1 \cap T_2} \subset \alpha_2 A_{T_2}$, thus  $x(\alpha_1 A_{T_1 \cap T_2}) \in C_1 \cap C_2$, therefore $C_1 \cap C_2\neq \emptyset$.

\bigskip\noindent
Suppose that $C_1 \cap C_2 \neq \emptyset$. Let $R$ be the smallest subset of $S$ such that $R_1 \cup R_2 \subset R \subset T_1 \cap T_2$ and  $\alpha_1^{-1} \alpha_2 \in A_R$. Set $C=C(\alpha_1 A_R, \alpha_1 A_{T_1 \cap T_2})$. It is easily seen that $C \subset C_1 \cap C_2$. Let $x(\beta A_X)$ be a vertex in $C_1 \cap C_2$. As before, we have $R_1 \cup R_2 \subset X \subset T_1 \cap T_2$, $\alpha_1^{-1} \alpha_2 \in A_X$, and $\alpha_1 A_X = \beta A_X$. By minimality of $R$, it follows that $R \subset X \subset T_1 \cap T_2$, thus $x(\beta A_X)$ is a vertex of $C$. So, $C_1 \cap C_2 \subset C$.

\bigskip\noindent
Suppose that $C_1,C_2$ span a cube. Let $C=C(\alpha A_R, \alpha A_T)$ be a cube containing $C_1$ and $C_2$. Then $\alpha A_R \subset \alpha_1 A_{R_1} \subset \alpha_1 A_{T_1} \subset \alpha A_T$ and $\alpha A_R \subset \alpha_2 A_{R_2} \subset \alpha_2 A_{T_2} \subset \alpha A_T$, thus $\alpha_1 A_{R_1} \cap \alpha_2 A_{R_2} \neq \emptyset$, because it contains $\alpha$. Moreover, $T_1 \cup T_2 \subset T$, thus $T_1 \cup T_2 \in \SS_{<\infty}$, since $T \in \SS_{<\infty}$.

\bigskip\noindent
Suppose that $\alpha_1 A_{R_1} \cap \alpha_2 A_{R_2} \neq \emptyset$ and $T_1 \cup T_2 \in \SS_{<\infty}$. Let $\alpha \in \alpha_1 A_{R_1} \cap \alpha_2 A_{R_2}$. Notice that $\alpha A_{R_1 \cap R_2} = \alpha_1 A_{R_1} \cap \alpha_2 A_{R_2}$. Then
\[\begin{array}{c}
\alpha A_{R_1 \cap R_2} \subset \alpha_1 A_{R_1} \subset \alpha_1 A_{T_1} \subset \alpha_1 A_{T_1 \cup T_2} = \alpha A_{T_1 \cup T_2}\,,\\
\noalign{\smallskip}
\alpha A_{R_1 \cap R_2} \subset \alpha_2 A_{R_2} \subset \alpha_2 A_{T_2} \subset \alpha_2 A_{T_1 \cup T_2} = \alpha A_{T_1 \cup T_2}\,,
\end{array}\]
thus $C(\alpha A_{R_1 \cap R_2}, \alpha A_{T_1 \cup T_2})$ contains $C_1$ and $C_2$. It is easily checked that, in this case, we have $C(\alpha A_{R_1 \cap R_2}, 
\alpha A_{T_1 \cup T_2}) = \Span (C_1,C_2)$.
\qed

\bigskip\noindent
{\bf Lemma 5.3.}
{\it There exists an algorithm which, given $X,Y,Z \in \SS_{<\infty}$, an element $\alpha \in A_Z$, and a word $\omega \in (\Sigma_Z \sqcup \Sigma_Z^{-1})^\ast$ which represents $\alpha$, decides whether $\alpha A_X \cap A_Y \neq \emptyset$ and, if yes, determines a word $\iota (\omega) \in (\Sigma_{Y \cap Z} \sqcup \Sigma_{Y \cap Z}^{-1})^\ast$ which represents an element in $\alpha A_X \cap A_Y$.}

\bigskip\noindent
{\bf Proof.}
Assume given $X,Y,Z \in \SS_{<\infty}$, $\alpha \in A_Z$, and $\omega \in (\Sigma_Z \sqcup \Sigma_Z^{-1})^\ast$ which represents $\alpha$. Set $w = \theta(\alpha)$. Write $w=w_0w_1$, where $w_0 \in W_Y$ and $w_1$ is $(Y,\emptyset)$-reduced. Recall that there is an algorithm which determines this decomposition. Note also that $w_0,w_1 \in W_Z$, since $w \in W_Z$ and $\lg(w) = \lg(w_0) + \lg(w_1)$.

\bigskip\noindent
First, we show that $wW_X \cap W_Y \neq \emptyset$ if and only if $w_0 \in w W_X$. Clearly, if $w_0 \in wW_X$, then $w_0 \in wW_X \cap W_Y$, thus $wW_X \cap W_Y \neq \emptyset$. Assume that $wW_X \cap W_Y \neq \emptyset$. Let $v \in wW_X \cap W_Y$. Let $v_0$ be the element of minimal length in $vW_X = wW_X$. Then $v_0$ is $(\emptyset,X)$-reduced and $v$ can be written $v=v_0v_1$, where $v_1 \in W_X$ and $\lg(v) = \lg(v_0) + \lg (v_1)$. Since $v \in W_Y$, this last equality implies that $v_0,v_1 \in W_Y$. On the other hand,  we can write $w=v_0w_2$, where $w_2 \in W_X$ and  $\lg(w) = \lg(v_0) + \lg(w_2)$, and we can write $w_2 = w_3w_4$, where $w_3 \in W_Y$ and $w_4$ is $(Y,\emptyset)$-reduced. So, $w=v_0 w_3 w_4$, $\lg(w) = \lg(v_0) + \lg (w_3) +\lg(w_4)$, $w_3 \in W_X \cap W_Y$, and $w_4$ is $(Y,\emptyset)$-reduced.  It follows that $w_0=v_0w_3$ and $w_1=w_4$, thus $w_0W_X = v_0W_X= wW_X$ (since $w_3 \in W_X$), therefore $w_0 \in wW_X \cap W_Y$.

\bigskip\noindent
We have $w_0 \in wW_X$ if and only if $w_0^{-1}w= w_1 \in W_X$, and, as said before, there exists an algorithm which decides whether $w_1 \in W_X$. Clearly, if $w W_X \cap W_Y = \emptyset$, then $\alpha A_X \cap A_Y = \emptyset$. So, we can (and do) assume that $w W_X \cap W_Y \neq \emptyset$. Then $w_1 \in W_X \cap W_Z$ and $w_0 \in w W_X \cap W_Y \cap W_Z$.

\bigskip\noindent
We compute reduced expressions $w_1=s_1 \cdots s_n$ and $w_0=t_1 \cdots t_m$, and we set $\tilde \tau(w_1) = \sigma_{s_1} \cdots \sigma_{s_n}$ and $\tilde \tau(w_0) = \sigma_{t_1} \cdots \sigma_{t_m}$. Note that $\tilde\tau(w_1) \in (\Sigma_{X \cap Z} \sqcup \Sigma_{X \cap Z}^{-1})^\ast$ and $\tilde\tau(w_0) \in (\Sigma_{Y \cap Z} \sqcup \Sigma_{Y \cap Z}^{-1})^\ast$ (since $w_1 \in W_{X \cap Z}$ and $w_0 \in W_{Y \cap Z}$). On the other hand, we set $\beta = \tau(w_0)^{-1} \alpha \tau(w_1)^{-1}$. Note that $\beta \in CA_Z$, $\beta A_X = \tau(w_0)^{-1} \alpha A_X$, and $\tau(w_0)^{-1}A_Y = A_Y$, thus $\alpha A_X \cap A_Y = \tau(w_0)(\beta A_X \cap A_Y)$.

\bigskip\noindent
Suppose that $\beta A_X \cap A_Y \neq \emptyset$. Let $\gamma \in \beta A_X \cap A_Y$. Set $v=\theta(\gamma)$. We have $v \in W_X \cap W_Y = W_{X \cap Y}$, thus $\tau(v) \in A_X \cap A_Y$, therefore $\gamma \tau(v)^{-1} \in \beta A_X \cap A_Y \cap CA = \beta CA_X \cap CA_Y$. So, $\beta A_X \cap A_Y \neq \emptyset$ if and only if $\beta CA_X \cap CA_Y \neq \emptyset$.

\bigskip\noindent
By the above, $\alpha A_X \cap A_Y \neq \emptyset$ if and only if $\beta CA_X \cap CA_Y \neq \emptyset$. By Lemma 4.6, $\beta CA_X \cap CA_Y \neq \emptyset$ if and only if $\pi_Y(\beta) = \pi_{Y\cap Z} (\beta) \in \beta CA_X$, that is, if and only if $\beta^{-1} \pi_{Y \cap Z} (\beta) \in CA_X$. Since $\beta^{-1} \pi_{Y \cap Z} (\beta) \in CA_Z$, we have $\beta^{-1} \pi_{Y \cap Z} (\beta) \in CA_X$ if and only if $\beta^{-1} \pi_{Y \cap Z} (\beta) \in CA_X \cap CA_Z = CA_{X \cap Z}$. Recall that there is an algorithm which, given  an element $\gamma \in CA_Z$ and a word $\mu \in (\Sigma_Z \sqcup \Sigma_Z^{-1})^\ast$ which represents $\gamma$, determines a word $\tilde \pi_{Y \cap Z} (\mu)$ which represents $\pi_{Y \cap Z} (\gamma)$ (this algorithm can be found just after the proof of Lemma 2.6). Set $\gamma = \beta^{-1} \pi_{Y \cap Z} (\beta)$ and $\mu = \tilde\tau(w_1) \omega^{-1} \tilde\tau(w_0) \cdot \tilde \pi_{Y \cap Z} (\tilde\tau(w_0)^{-1} \omega \tilde\tau(w_1)^{-1})$. Then $\mu$ is a word in $(\Sigma_Z \sqcup \Sigma_Z^{-1})^\ast$ which represents $\gamma$, and we have $\gamma \in CA_{X \cap Z}$ if and only if $\pi_{X \cap Z} (\gamma) = \gamma$. Finally, we can check whether $\pi_{X \cap Z} (\gamma) = \gamma$ applying to $\mu$ and $\tilde \pi_{X \cap Z} (\mu)$ a solution to the word problem in $A_Z$.

\bigskip\noindent
Suppose that $\alpha A_X \cap A_Y \neq \emptyset$. Then $\pi_{Y \cap Z} (\beta) \in \beta A_X \cap A_Y$, thus $\tau(w_0) \pi_{Y \cap Z} (\beta) \in \alpha A_X \cap A_Y$, and $\iota(\omega) = \tilde\tau(w_0) \tilde \pi_{Y \cap Z}(\tilde\tau(w_0)^{-1} \omega \tilde\tau(w_1)^{-1})$ is a word in $(\Sigma_{Y \cap Z} \sqcup \Sigma_{Y \cap Z}^{-1})^\ast$ which represents $\tau(w_0) \pi_{Y \cap Z} (\beta)$.
\qed

\bigskip\noindent
{\bf Lemma 5.4.}
{\it There exists an algorithm which, given two cubes $C_1=C(\alpha_1 A_{R_1}, \alpha_1 A_{T_1})$ and $C_2=C(\alpha_2 A_{R_2}, \alpha_2 A_{T_2})$, a word $\omega_1 \in (\Sigma \sqcup \Sigma^{-1})^\ast$ which represents $\alpha_1$, an element $X \in \SS_{<\infty}$, and a word $\nu \in (\Sigma_X \sqcup \Sigma_X^{-1})^\ast$ such that $\omega_2 = \omega_1 \nu$ represents $\alpha_2$,
\begin{itemize}
\item
decides whether $C_1 \cap C_2 \neq \emptyset$ and, if yes, determines $R \in \SS_{<\infty}$ such that $C_1 \cap C_2 = C(\alpha_1 A_R, \alpha_1 A_{T_1\cap T_2})$,
\item
decides whether $C_1$ and $C_2$ span a cube, and, if yes, determines a word $\mu \in (\Sigma_{R_1\cap X} \sqcup \Sigma_{R_1\cap X}^{-1})^\ast$ such that $\Span (C_1,C_2) = C(\alpha A_{R_1 \cap R_2}, \alpha A_{T_1 \cup T_2})$, where $\alpha$ is the element of $A$ represented by $\omega_1 \mu$.
\end{itemize}}

\bigskip\noindent
{\bf Proof.} 
Recall that $C_1 \cap C_2 \neq \emptyset$ if and only if $R_1 \cup R_2 \subset T_1 \cap T_2$ and $\alpha_1^{-1} \alpha_2 \in A_{T_1\cap T_2}$ (see Lemma 5.2). It is easy to check whether $R_1 \cup R_2 \subset T_1 \cap T_2$ or not. Set $\beta = \alpha_1^{-1} \alpha_2$. Then $\beta$ is represented by $\nu \in (\Sigma_X \sqcup \Sigma_X^{-1})^\ast$. In particular, $\beta \in A_X$, thus we have $\beta \in A_{T_1 \cap T_2}$ if an only if $\beta \in A_{T_1 \cap T_2} \cap A_X = A_{T_1 \cap T_2 \cap X}$. Recall that $A_X$ has a solution to the word problem (since $X \in \SS_{<\infty}$), thus, by Proposition 2.7, there is an algorithm which decides whether $\beta \in A_{T_1 \cap T_2 \cap X}$. So, there exists an algorithm which decides whether $C_1 \cap C_2 \neq \emptyset$.

\bigskip\noindent
Suppose that $C_1 \cap C_2 \neq \emptyset$. Let $R \in \SS_{<\infty}$ such that $R_1 \cup R_2 \subset R \subset T_1 \cap T_2$. We have $\beta \in A_R$ if and only if $\beta \in A_R \cap A_X = A_{R\cap X}$, and, by Proposition 2.7, there is an algorithm which decides whether $\beta \in A_{R \cap X}$. So, there exists an algorithm which determines the smallest $R \in \SS_{<\infty}$ such that $R_1 \cup R_2 \subset R \subset T_1 \cap T_2$ and $\beta \in A_R$. By Lemma 5.2, this $R$ satisfies $C_1 \cap C_2 = C(\alpha_1 A_R, \alpha_1 A_{T_1 \cap T_2})$.

\bigskip\noindent
By Lemma 5.2, $C_1$ and $C_2$ span a cube if and only if $\alpha_1 A_{R_1} \cap \alpha_2 A_{R_2} \neq \emptyset$ and $T_1 \cup T_2 \in \SS_{<\infty}$. It is easy to check whether $T_1 \cup T_2 \in \SS_{<\infty}$ or not. We have $\alpha_1 A_{R_1} \cap \alpha_2 A_{R_2} \neq \emptyset$ if and only if 
\break
$A_{R_1} \cap \beta A_{R_2} \neq \emptyset$. By Lemma 5.3, there is an algorithm which decides whether $A_{R_1} \cap \beta A_{R_2} \neq \emptyset$. Suppose $A_{R_1} \cap \beta A_{R_2} \neq \emptyset$. Then this algorithm also determines a word $\mu \in (\Sigma_{R_1\cap X} \sqcup \Sigma_{R_1\cap X}^{-1})^\ast$ which represents an element $\gamma \in A_{R_1} \cap \beta A_{R_2}$. Set $\alpha = \alpha_1 \gamma$. Then $\omega_1 \mu$ represents $\alpha \in \alpha_1 A_{R_1} \cap \alpha_2 A_{R_2}$, and $\Span(C_1,C_2) = C(\alpha A_{R_1 \cap R_2}, \alpha A_{T_1 \cup T_2})$.
\qed

\bigskip\noindent
Let $n \ge 1$. A {\it cube prepath of length $n\ge 1$} is defined to be a sequence of $n$ words $\omega_1, \nu_2, \dots, \nu_n \in (\Sigma \sqcup \Sigma^{-1})^\ast$ together with $2n+2$ elements of $\SS_{<\infty}$, $R_1, \dots, R_n, T_1, \dots, T_n, X, Y \in \SS_{<\infty}$, such that $R_i \subset T_i$ for all $1 \le i \le n$, $R_{i-1} \cup R_i \subset T_{i-1} \cap T_i$ and $\nu_i \in (\Sigma_{T_{i-1} \cap T_i} \sqcup \Sigma_{T_{i-1} \cap T_i}^{-1})^\ast$ for all $2 \le i \le n$, $R_1 \subset X \subset T_1$, and $R_n \subset Y \subset T_n$. Such a cube prepath is denoted by
\[
\PP = ((\omega_1, \nu_2, \dots, \nu_n),(R_1, \dots, R_n),(T_1, \dots, T_n),(X,Y))\,.
\]
A {\it cube prepath of length $0$} is a word $\omega \in (\Sigma \sqcup \Sigma^{-1})^\ast$ together with an element $X \in \SS_{<\infty}$. This cube prepath is denoted by $\PP=((\omega),(),(),(X,X))$.

\bigskip\noindent
For $1 \le i \le n$, we denote by $\alpha_i$ the element of $A$ represented by $\omega_1 \nu_2 \cdots \nu_i$, and we set $C_i = C(\alpha_i A_{R_i}, \alpha_i A_{T_i})$. We also set $x = x(\alpha_1 A_X)$ and $y=x(\alpha_n A_Y)$. By Lemma 5.2, we have $C_i \cap C_{i+1} \neq \emptyset$ for all $1 \le i \le n-1$, thus $\CC=(C_1, \dots, C_n)$ is a cube path in $\Phi$ from $x$ to $y$. We say that $\CC$ is the {\it geometric realization} of $\PP$, and that $\PP$ is a cube prepath from $x$ to $y$. We say that $\PP$ is {\it normal} if $\CC$ is normal. If $\PP$ is of length $0$, then $x=y$ and $\CC = ()$ is the constant cube path on $x$.

\bigskip\noindent
The following theorem is the main result of this section. Our solution to the word problem in $A$ will be a straightforward consequence of it. 

\bigskip\noindent
{\bf Theorem 5.5.}
{\it There exists an algorithm which, given two vertices $x$ and $y$ in $\Phi$, and a cube prepath $\PP$ from $x$ to $y$, determines a normal cube prepath $\PP'$ from $x$ to $y$.}

\bigskip\noindent
{\bf Proof.}
For $\omega \in (\Sigma \sqcup \Sigma^{-1})^\ast$, we denote by $\bar \omega$ the element of $A$ represented by $\omega$. We assume given a cube prepath 
\[
\PP = ((\omega_1, \nu_2, \dots, \nu_n), (R_1, \dots, R_n), (T_1, \dots, T_n),(X,Y))\,.
\]
We set $\alpha_i = \bar \omega_1 \bar \nu_2 \cdots \bar \nu_i$ and $C_i = C(\alpha_i A_{R_i}, \alpha_i A_{T_i})$ for all $1 \le i \le n$. We also set $x = x(\alpha_1 A_X)$ and $y = x(\alpha_n A_Y)$. So, $\PP$ is a cube prepath from $x$ to $y$. Our aim is to make an effective construction of a normal cube prepath
\[
\PP' = ((\omega_1', \nu_2', \dots, \nu_m'), (R_1', \dots, R_m'), (T_1', \dots, T_m') ,(X,Y))
\]
from $x$ to $y$ such that $m \le n$. Set $\alpha_m' = \bar \omega_1' \bar \nu_2' \cdots \bar \nu_m'$. Then this construction also provides a word $\mu' \in (\Sigma_Y \sqcup \Sigma_Y^{-1})^\ast$ which represents $( \alpha_m')^{-1} \alpha_n$ (this word will be used in the inductive step). We argue by double induction, on $n = \lg (\PP)$, and on $\dim (C_n)$. The case $n=0$ being trivial, we can and do assume $n \ge 1$.

\bigskip\noindent
Suppose $n=1$. Set $C_1' = \Span(x,y) = C(\alpha_1 A_{X \cap Y}, \alpha_1 A_{X \cup Y})$. If $X = Y$, then $x=y$ and $\PP'= ((\omega_1),(),(),(X,X))$ is a normal cube prepath from $x$ to $y=x$. In this case we set $\mu'=1$. Suppose that $X \neq Y$. Then $\CC'=(C_1')$ is a normal cube path  from $x$ to $y$. Set $\PP'= ((\omega_1),(X \cap Y), (X \cup Y),(X,Y))$. Then $\PP'$ is a (normal) cube prepath whose geometric realization is $\CC'$. We also set $\mu'=1$ in this case.

\bigskip\noindent
Suppose $n \ge 2$ and $\dim (C_n) = 0$, plus the inductive hypothesis (on $n$). Then $R_n = Y = T_n$ and $y =x(\alpha_{n-1} A_Y)$. We also have $\nu_n \in (\Sigma_{T_{n-1} \cap T_n} \sqcup \Sigma_{T_{n-1} \cap T_n}^{-1})^\ast = (\Sigma_Y \sqcup \Sigma_Y^{-1})^\ast$. Set
\[
\PP^{(1)} = (( \omega_1, \nu_2, \dots, \nu_{n-1}), (R_1, \dots, R_{n-1}), (T_1, \dots, T_{n-1}),(X,Y))\,.
\]
Then $\PP^{(1)}$ is a cube prepath from $x$ to $y$. By induction, we can construct a normal cube prepath
\[
\PP' = (( \omega_1', \nu_2', \dots, \nu_m'),(R_1', \dots, R_m'),(T_1', \dots, T_m'),(X,Y))
\]
from $x$ to $y$ such that $m \le n-1$. Set $\alpha_m' = \bar \omega_1' \bar \nu_2' \cdots \bar \nu_m'$. This construction also provides a word $\mu^{(1)} \in (\Sigma_Y \sqcup \Sigma_Y^{-1})^\ast$ which represents $(\alpha_m')^{-1} \alpha_{n-1}$. Set $\mu' = \mu^{(1)} \nu_n$. Then $\mu' \in (\Sigma_Y \sqcup \Sigma_Y^{-1})^\ast$, and $\mu'$ represents $(\alpha_m')^{-1} \alpha_n$.

\bigskip\noindent
Now, we assume that $n \ge 2$, $\dim (C_n) \ge 1$, plus the inductive hypothesis. The remainder of the construction is divided into 4 steps.
\begin{itemize}
\item
In Step 1 we prove that our study can be reduced to the case where (1) $\dim (C_i) \ge 1$ for all $1 \le  i\le n-1$, (2) $C_i \cap C_{i+1}$ is a vertex, denoted by $x_i$, for all $1 \le i \le n-2$, (3) $C_i=\Span(x_{i-1},x_i)$ for all $1 \le i\le n-1$, where $x_0=x$ and $x_{n-1}\in C_{n-1} \cap C_n$, (4) $\Star(C_i) \cap C_{i+1} = \{x_i\}$ for all $1 \le i \le n-2$.
\item
In Step 2 we prove that our study can be reduced to the case  where  $C_n = \Span(x_{n-1},y)$ (plus the conditions of Step 1).
\item
In Step 3 we consider a vertex $z'$ of $C_n$ such that $\Span(x_{n-2},z')$ is nonempty and of minimal dimension, and we prove that we can reduce our study to the case where $C_{n-1} = \Span(x_{n-2},z')$ and $C_n=\Span(z',y)$ (plus the conditions of Steps 1 and 2).
\item
In Step 4 we prove that, under the conditions of the previous steps, we have $\Star(C_{n-1}) \cap C_n = \{z'\}$.
\end{itemize}

\bigskip\noindent
{\bf Step 1.}
Choose a vertex $z$ in $C_{n-1} \cap C_n$. Let $Z$ be the element of $\SS_{<\infty}$ such that $R_{n-1} \cup R_n \subset Z \subset T_{n-1} \cap T_n$ and $z = x(\alpha_{n-1} A_Z) = x(\alpha_n A_Z)$. Set 
\[
\PP^{(1)} = ((\omega_1, \nu_2, \dots, \nu_{n-1}), (R_1, \dots, R_{n-1}), (T_1, \dots, T_{n-1}),(X,Z))\,.
\]
Then $\PP^{(1)}$ is a cube prepath of length $n-1$ from $x$ to $z$. By induction, we can construct a normal cube prepath 
\[
\PP^{(2)} = ((\omega_1^{(2)}, \nu_2^{(2)}, \dots, \nu_{m-1}^{(2)}), (R_1^{(2)}, \dots, R_{m-1}^{(2)}), (T_1^{(2)}, \dots, T_{m-1}^{(2)}), (X,Z))
\]
from $x$ to $z$ such that $m-1 \le n-1$. Set $\alpha_{m-1}^{(2)} = \bar \omega_1^{(2)} \bar \nu_2^{(2)} \cdots \bar \nu_{m-1}^{(2)}$. Then this construction also provides a word $\mu^{(2)} \in (\Sigma_Z \sqcup \Sigma_Z^{-1})^\ast$ which represents $(\alpha_{m-1}^{(2)})^{-1} \alpha_{n-1}$. The word $\mu^{(2)} \nu_n$ represents $(\alpha_{m-1}^{(2)})^{-1} \alpha_n$ and belongs to $(\Sigma_{T_{n-1} \cap T_n} \sqcup \Sigma_{T_{n-1} \cap T_n}^{-1})^\ast$. Set $C_{m-1}^{(2)} = C(\alpha_{m-1}^{(2)} A_{R_{m-1}^{(2)}}, \alpha_{m-1}^{(2)} A_{T_{m-1}^{(2)}})$. Since $C_{m-1}^{(2)} \cap C_n \neq \emptyset$, by Lemma 5.2, we have $(\alpha_{m-1}^{(2)})^{-1} \alpha_n \in A_{T_{m-1}^{(2)} \cap T_n}$, thus $(\alpha_{m-1}^{(2)})^{-1} \alpha _n \in A_{T_{m-1}^{(2)} \cap T_{n-1} \cap T_n}$. By Proposition 2.7, there is an algorithm which determines a word $\nu_m^{(2)}= \kappa (\mu^{(2)} \nu_n) \in (\Sigma_{T_{m-1}^{(2)} \cap T_{n-1} \cap T_n} \sqcup \Sigma_{T_{m-1}^{(2)} \cap T_{n-1} \cap T_n}^{-1})^\ast$ which represents $(\alpha_{m-1}^{(2)})^{-1} \alpha_n$. Set
\[
\PP^{(3)} = (( \omega_1^{(2)}, \nu_2^{(2)}, \dots, \nu_{m-1}^{(2)},\nu_m^{(2)}),(R_1^{(2)}, \dots, R_{m-1}^{(2)},R_n), (T_1^{(2)}, \dots, T_{m-1}^{(2)}, T_n), (X,Y))\,.
\]
Then $\PP^{(3)}$ is a cube prepath from $x$ to $y$. Moreover, $\bar \omega_1^{(2)} \bar \nu_2^{(2)} \cdots \bar \nu_{m-1}^{(2)} \bar \nu_m^{(2)} = \alpha_n$. 

\bigskip\noindent
By induction, we can suppose that $m=n$, and, up to replacing $\PP$ by $\PP^{(3)}$, we can assume that:
\begin{itemize}
\item
$\dim (C_i) \ge 1$ for all $1 \le  i\le n-1$,
\item
$C_i \cap C_{i+1}$ is a vertex, denoted by $x_i$, for all $1 \le i \le n-2$,
\item
$C_i=\Span(x_{i-1},x_i)$ for all $1 \le i\le n-1$, where $x_0=x$ and $x_{n-1}=z$,
\item
$\Star(C_i) \cap C_{i+1} = \{x_i\}$ for all $1 \le i \le n-2$.
\end{itemize}

\bigskip\noindent
{\bf Step 2.} 
Set $C_n^{(1)} = \Span (z,y) = C(\alpha_n A_{Z \cap Y}, \alpha_n A_{Z \cup Y})$. Since $C_{n-1} \cap C_n^{(1)} \neq \emptyset$, we have $\alpha_{n-1}^{-1} \alpha_n \in A_{T_{n-1} \cap (Z \cup Y)}$. Moreover, $T_{n-1} \cap (Z \cup Y) \subset T_{n-1} \cap T_n$, because $Z \cup  Y \subset T_n$. By Proposition 2.7, there is an algorithm which determines a word $\nu_n^{(1)} = \kappa (\nu_n) \in (\Sigma_{T_{n-1} \cap (Z \cup Y)} \sqcup \Sigma_{T_{n-1} \cap (Z \cup Y)}^{-1})^\ast$ which represents $\alpha_{n-1}^{-1} \alpha_n$. Set 
\[
\PP^{(1)} = ((\omega_1, \nu_2, \dots, \nu_{n-1}, \nu_n^{(1)}), (R_1, \dots, R_{n-1}, Z \cap Y), (T_1, \dots, T_{n-1}, Z \cup Y), (X,Y))\,.
\]
Then $\PP^{(1)}$ is a cube prepath from $x$ to $y$. Furthermore, $\alpha_n = \bar \omega_1 \bar \nu_2 \cdots \bar \nu_{n-1} \bar \nu_n^{(1)}$. So, up to replacing $\PP$ by $\PP^{(1)}$, we can assume that $C_n= \Span(z,y)$.

\bigskip\noindent
{\bf Step 3.}
Consider the set of vertices $x (\alpha_n A_U)$ of $C_n$ such that $\Span(x_{n-2},x(\alpha_n A_U))$ is nonempty. In other words, consider the set $\UU = \{ U \in \SS_{<\infty}; R_n \subset U \subset T_n \text{ and } \Span(x_{n-2},x(\alpha_n A_U)) \neq \emptyset\}$. This set is nonempty since it contains $Z$. Choose $Z' \in \UU$ such that the dimension of $\Span(x(\alpha_n A_{Z'}), y) \subset C_n$ is minimal, and set $z'=x(\alpha_n A_{Z'})$, $C_{n-1}^{(1)} = \Span(x_{n-2},z')$, $C_n^{(1)}= \Span(z',y)$, and 
\[
\CC^{(1)}=(C_1,\dots, C_{n-2}, C_{n-1}^{(1)}, C_n^{(1)})\,.
\]
Note that $\CC^{(1)}$ is a cube path from $x$ to $y$. Note also that, thanks to Lemma 5.4, it is easy to determine in an algorithmic way the set $\UU$ as well as the element $Z' \in \UU$.

\bigskip\noindent
Since $C_{n-2} \cap C_{n-1} = \{x_{n-2}\}$, by Lemma 5.2, we have $R_{n-2} \cup R_{n-1} = T_{n-2} \cap T_{n-1}$ and $x_{n-2} = x(\alpha_{n-2} A_{Z_{n-2}}) = x(\alpha_{n-1} A_{Z_{n-2}})$, where $Z_{n-2} = R_{n-2} \cup R_{n-1} = T_{n-2} \cap T_{n-1}$. Set $R_{n-1}^{(1)} = Z_{n-2} \cap Z'$ and $T_{n-1}^{(1)} = Z_{n-2} \cup Z'$. Recall also that $\nu_n \in (\Sigma_{T_{n-1} \cap T_n} \sqcup \Sigma_{T_{n-1} \cap T_n}^{-1})^\ast$. By Lemma~5.4, there is an algorithm which determines a word $\mu^{(1)} \in (\Sigma_{T_{n-1} \cap T_n \cap Z_{n-2}} \sqcup \Sigma_{T_{n-1} \cap T_n \cap Z_{n-2}}^{-1})^\ast$ such that $C_{n-1}^{(1)} = C(\alpha_{n-1}^{(1)} A_{R_{n-1}^{(1)}}, \alpha_{n-1}^{(1)} A_{T_{n-1}^{(1)}})$, where $\alpha_{n-1}^{(1)}$ is the element of $A$ represented by 
\break
$\omega_1 \nu_2 \cdots \nu_{n-2} \nu_{n-1} \mu^{(1)}$. Observe that $\nu_{n-1} \mu^{(1)}$ lies in $(\Sigma_{T_{n-2} \cap T_{n-1}} \sqcup \Sigma_{T_{n-2} \cap T_{n-1}}^{-1})^\ast$ and represents $\alpha_{n-2}^{-1} \alpha_{n-1}^{(1)}$. On the other hand, by Lemma 5.2, $\alpha_{n-2}^{-1} \alpha_{n-1}^{(1)} \in A_{T_{n-2} \cap T_{n-1}^{(1)}}$, thus $\alpha_{n-2}^{-1} \alpha_{n-1}^{(1)} \in A_{T_{n-2} \cap T_{n-1} \cap T_{n-1}^{(1)}}$. By Proposition 2.7, there is an algorithm which determines a word $\nu_{n-1}^{(1)} = \kappa (\nu_{n-1} \mu^{(1)}) \in (\Sigma_{T_{n-2} \cap T_{n-1} \cap T_{n-1}^{(1)}} \sqcup \Sigma_{T_{n-2} \cap T_{n-1} \cap T_{n-1}^{(1)}}^{-1})^\ast$ which represents $\alpha_{n-2}^{-1} \alpha_{n-1}^{(1)}$.

\bigskip\noindent
Recall that $y=x(\alpha_n A_Y)$. Set $R_n^{(1)} = Z' \cap Y$ and $T_n^{(1)} = Z' \cup Y$. We have $C_n^{(1)} = C(\alpha_n A_{R_n^{(1)}}, 
\break
\alpha_n A_{T_n^{(1)}})$. The word $(\mu^{(1)})^{-1} \nu_n$ lies in $(\Sigma_{T_{n-1} \cap T_n} \sqcup \Sigma_{T_{n-1} \cap T_n}^{-1})^\ast$ and represents $(\alpha_{n-1}^{(1)})^{-1} \alpha_n$. On the other hand, by Lemma 5.2, we have $(\alpha_{n-1}^{(1)})^{-1} \alpha_n \in A_{T_{n-1}^{(1)} \cap T_n^{(1)}}$, thus $(\alpha_{n-1}^{(1)})^{-1} \alpha_n \in A_{T_{n-1} \cap T_n \cap T_{n-1}^{(1)} \cap T_n^{(1)}}$. By Proposition 2.7, there is an algorithm which determines a word $\nu_n^{(1)} = \kappa ((\mu^{(1)})^{-1} \nu_n) \in (\Sigma_{T_{n-1} \cap T_n \cap T_{n-1}^{(1)} \cap T_n^{(1)}} \sqcup \Sigma_{T_{n-1} \cap T_n \cap T_{n-1}^{(1)} \cap T_n^{(1)}}^{-1})^\ast$ which represents $(\alpha_{n-1}^{(1)})^{-1} \alpha_n$.

\bigskip\noindent
Set
\begin{multline*}
\PP^{(1)} = ((\omega_1, \nu_2, \dots, \nu_{n-2}, \nu_{n-1}^{(1)}, \nu_n^{(1)}), (R_1, \dots, R_{n-2}, R_{n-1}^{(1)}, R_n^{(1)}), \\
(T_1, \dots, T_{n-2}, T_{n-1}^{(1)}, T_n^{(1)}), (X,Y))\,.
\end{multline*}
Then $\PP^{(1)}$ is a cube prepath whose geometric realization is $\CC^{(1)}$. Moreover, we have $\alpha_n = \bar \omega_1 \bar \nu_2 \cdots \bar \nu_{n-2} \bar \nu_{n-1}^{(1)} \bar \nu_n^{(1)}$.

\bigskip\noindent
If $z' \neq z$, then $C_n^{(1)} = \Span(z',y)$ is a proper face of $C_n=\Span(z,y)$. In that case, we can conclude using the induction on $\dim (C_n^{(1)})$. So, we can and do assume that $z'=z$. Then $C_n = \Span(z,y) = C_n^{(1)}$, and $C_{n-1} = \Span(x_{n-2},z) = C_{n-1}^{(1)}$.

\bigskip\noindent
{\bf Step 4.}
Now, we show that, under the hypothesis deduced from the previous steps, $\PP$ is a normal cube prepath. By Step 1, we have
\begin{itemize}
\item
$\dim (C_i) \ge 1$ for all $1 \le  i\le n-1$,
\item
$C_i \cap C_{i+1}$ is a vertex, denoted by $x_i$, for all $1 \le i \le n-2$,
\item
$C_i=\Span(x_{i-1},x_i)$ for all $1 \le i\le n-1$, where $x_0=x$ and $x_{n-1}=z$,
\item
$\Star(C_i) \cap C_{i+1} = \{x_i\}$ for all $1 \le i \le n-2$.
\end{itemize}
We also have $C_n = \Span(z,y)$ (by Step 2) and $\dim (C_n) \ge 1$ (by induction). It remains to show that $\Star(C_{n-1}) \cap C_n = \{z\}$.

\bigskip\noindent
Recall that $\UU = \{ U \in \SS_{<\infty}; R_n \subset U \subset T_n \text{ and } \Span(x_{n-2}, x(\alpha_n A_U)) \neq \emptyset\}$. Recall also that, according to Step 3, $Z=Z'$ is such that the dimension of $\Span(x(\alpha_n A_Z),y)$ is minimal among the dimensions of $\Span(x(\alpha_n A_U),y)$, with $U \in \UU$. If $U \in \UU$, then $R_n \subset U \subset T_n$, thus
\[
R_n = Z \cap Y \subset U \cap Y \subset U \cup Y \subset T_n= Z \cup Y\,.
\]
By minimality of the dimension of $\Span(x(\alpha_n A_Z),y)$, it follows that $Z \cap Y = U \cap Y$ and $Z \cup Y = U \cup Y$, thus $U=Z$. So, $\UU = \{Z\}$. Now, if $x(\alpha_n A_U)$ is a vertex in $\Star(C_{n-1}) \cap C_n$, then $U \in \UU = \{ Z\}$, thus $x(\alpha_n A_U)=z$.
\qed

\bigskip\noindent
{\bf Theorem 5.6.}
{\it The group $A$ has a solution to the word problem.}

\bigskip\noindent
{\bf Proof.}
Let $\omega = \sigma_{s_1}^{\varepsilon_1} \cdots \sigma_{s_n}^{\varepsilon_n}$ be a word in $(\Sigma \sqcup \Sigma^{-1})^\ast$. Set $\alpha = \bar \omega$. We would like to determine whether $\alpha =1$ or not. Set $x = x(1 A_\emptyset)$ and $y = x(\alpha A_\emptyset)$. We have $\alpha =1$ if and only if $x=y$. On the other hand, we have $x=y$ if and only if the unique normal cube path in $\Phi$ from $x$ to $y$ is of length $0$.

\bigskip\noindent
Set $\omega_1 = 1$, $X=\emptyset$, and $Y = \emptyset$. For $1 \le i \le n$, set $\nu_{2i} = \sigma_{s_i}^{\varepsilon_i}$, $R_{2i}=\emptyset$, and $T_{2i} = \{s_i\}$. For $2 \le i \le n$, set $\nu_{2i-1} = 1$, $R_{2i-1} = \emptyset$, and $T_{2i-1} = \{s_i\}$. Set $R_1 = \emptyset$, and $T_1 = \{s_1\}$ (there is no $\nu_1$). Then
\[
\PP = ((\omega_1, \nu_2, \dots, \nu_{2n}), (R_1, \dots, R_{2n}), (T_1, \dots, T_{2n}), (X,Y))
\]
is a cube prepath from $x$ to $y$. By Theorem 5.5, there is an algorithm which determines from $\PP$ a normal cube prepath $\PP'$ from $x$ to $y$. Then $\alpha =1$ if and only if $\PP'$ is of length $0$.
\qed




\section{Virtual braid groups}

The {\it virtual braid group} on $n$ strands, denoted $VB_n$, is defined by the presentation with generators $\sigma_1, \dots, \sigma_{n-1}, \tau_1, \dots, \tau_{n-1}$, and relations
\[ \begin{array}{cl}
\tau_i^2=1 &\ \text{if } 1 \le i \le n-1\\
\sigma_i \sigma_j = \sigma_j \sigma_i\,,\ 
\sigma_i \tau_j = \tau_j \sigma_i\,,  \text{ and }
\tau_i \tau_j = \tau_j \tau_i
&\ \text{if } |i-j| \ge 2\\
\sigma_i \sigma_j \sigma_i = \sigma_j \sigma_i \sigma_j\,,\
\sigma_i \tau_j \tau_i = \tau_j \tau_i \sigma_j\,,\text{ and }
\tau_i \tau_j \tau_i = \tau_j \tau_i \tau_j
& \ \text{if } |i-j|=1
\end{array}\]
Virtual braids can be viewed as braid diagrams that admit virtual crossings of the same type as the virtual crossings in the virtual links. They were introduced by Kauffman in \cite{Kauff1} at the same time as the virtual links. Like for the classical links and braids, every virtual link is the closure of a virtual braid, and the closures of two virtual braids are equivalent if and only if the braids are related by a finite sequence of moves called virtual Markov moves (see \cite{Kamad1} and \cite{KauLam1}).

\bigskip\noindent
Let $\SSS_n$ denote the symmetric group of $\{1, \dots, n\}$. There is an epimorphism $\theta : VB_n \to \SSS_n$ which sends $\sigma_i$ to $1$ and $\tau_i$ to $(i,i+1)$ for all $1 \le i \le n-1$. This epimorphism admits a section $\iota : \SSS_n \to VB_n$ which sends $(i,i+1)$ to $\tau_i$ for all $1 \le i \le n-1$. Let $K_n$ denote the kernel of $\theta : VB_n \to \SSS_n$. Then, by the above, $VB_n$ is a semidirect product $K_n \rtimes \SSS_n$. A precise description of this semidirect product was determined by Rabenda in his {\it Mémoire de DEA} (Master degree thesis) \cite{Raben1} (see also \cite{BarBel1}).

\bigskip\noindent
{\bf Remark.}
There is another ``natural'' epimorphism $\theta' : VB_n \to \SSS_n$ which sends $\sigma_i$ and $\tau_i$ to $(i,i+1)$ for all $1 \le i \le n-1$. This epimorphism also admits a section $\iota' : \SSS_n \to VB_n$ which sends $(i,i+1)$ to $\tau_i$ for all $1 \le i \le n-1$. The kernel of $\theta'$ is called the {\it virtual pure braid group}, and is denoted by $VP_n$. It is different from the group $K_n$ studied in the present section (see \cite{BarBel1}).

\bigskip\noindent
{\bf Proposition 6.1} (Rabenda \cite{Raben1}).
{\it For $1 \le i < j \le n$, we set
\[ \begin{array}{rcl}
\delta_{i,j} &=& \tau_i \tau_{i+1} \cdots \tau_{j-2} \sigma_{j-1} \tau_{j-2} \cdots \tau_{i+1} \tau_i\\
\noalign{\smallskip}
\delta_{j,i} &=& \tau_i \tau_{i+1} \cdots \tau_{j-2}\tau_{j-1} \sigma_{j-1} \tau_{j-1} \tau_{j-2} \cdots \tau_{i+1} \tau_i\\
\end{array}\]
Then $K_n$ has a presentation with generators $\delta_{i,j}$, $1 \le i \neq j \le n$, and relations
\[\begin{array}{cl}
\delta_{i,j} \delta_{k,l} = \delta_{k,l} \delta_{i,j} &\quad \text{for } i,j,k,l \text{ distincts}\\
\noalign{\smallskip}
\delta_{i,j} \delta_{j,k} \delta_{i,j} = \delta_{j,k} \delta_{i,j} \delta_{j,k} &\quad \text{for } i,j,k \text{ distincts}
\end{array}\]
Moreover, the action of $\SSS_n$ on $K_n$ is given by
\[
w \cdot \delta_{i,j} = \delta_{w(i), w(j)}\text{ for }1 \le i \neq j \le n \text{ and } w \in \SSS_n\,.
\]}

\bigskip\noindent
Let $\Gamma_{VB,n}$ be the Coxeter graph defined by the following data. 
\begin{itemize}
\item
The set of vertices of $\Gamma_{VB,n}$ is $S=\{x_{i,j} ; 1 \le i \neq j \le n\}$.
\item
If $i,j,k,l$ are distinct, then $x_{i,j}$ and $x_{k,l}$ are not joined by any edge.
\item
If $i,j,k$ are distinct, then $x_{i,j}$ and $x_{j,k}$ are joined an ordinary edge (labeled by $3$), $x_{i,j}$ and $x_{i,k}$ are joined by an edge labeled by $\infty$, and $x_{i,k}$ and $x_{j,k}$ are joined by an edge labeled by $\infty$.
\item
If $i,j$ are distinct, then $x_{i,j}$ and $x_{j,i}$ are joined by an edge labeled by $\infty$.
\end{itemize}
By Proposition 6.1, $K_n$ is an Artin-Tits group of $\Gamma_{VB,n}$. 

\bigskip\noindent
{\bf Proposition 6.2.}
{\it Let $\SS_{<\infty}$ denote the set of $X \subset S$ such that $(\Gamma_{VB,n})_X$ is free of infinity. Then $(\Gamma_{VB,n})_X$ is of type $K(\pi,1)$ and $(K_n)_X$ has a solution to the word problem for all $X \in \SS_{<\infty}$.}

\bigskip\noindent
{\bf Proof.}
Let $X \in \SS_{<\infty}$. It is easily seen that every edge of $(\Gamma_{VB,n})_X$ is  ordinary ({\it i.e.} not labeled), and the valence of every vertex is less or equal to $2$. It follows that $(\Gamma_{VB,n})_X$ is the disjoint union of Coxeter graphs of type $A$ and $\tilde A$. By \cite{FadNeu1}, any Coxeter graph of type $A$ is of type $K(\pi,1)$ (see also \cite{Delig1}), and, by \cite{Okone1}, any Coxeter graph of type $\tilde A$ is of type $K(\pi,1)$ (see also \cite{ChaPei1}). The Artin-Tits groups of type $A$ are the Artin braid groups, and these are known to have several (fast) solutions to the word problem (see \cite{Artin2}, \cite{Artin1}, \cite{Garsi1}, for instance). On the other hand, there are natural embeddings of the Artin-Tits groups of type $\tilde A$ into the Artin-Tits groups of type $B$ (and, therefore, into the Artin braid groups) (see \cite{Allco1}, \cite{Dieck1}, \cite{KenPei1}). One can easily solve the word problem in the Artin-Tits groups of type $\tilde A$ using these embeddings. Another solution to the word problem in the Artin-Tits groups of type $\tilde A$ can be found in \cite{Digne1}.
\qed

\bigskip\noindent
{\bf Corollary 6.3.}
{\it \begin{enumerate}
\item
$K_n$ and $VB_n$ have solutions to the word problem.
\item
$K_n$ has a classifying space of dimension $n-1$, and its cohomological dimension is $n-1$. In particular, $K_n$ is of FP type and is torsion free.
\item
$VB_n$ is virtually torsion free, and its virtual cohomological dimension is $n-1$.
\end{enumerate}}

\bigskip\noindent
{\bf Proof.}
The fact that $K_n$ has a solution to the word problem is a straightforward consequence of Theorem 5.6 and Proposition 6.2. Since $K_n$ is a finite index subgroup of $VB_n$, it follows that $VB_n$ has also a solution to the word problem.

\bigskip\noindent
Let $d= \max\{ |X|; X \in \SS^f\}$. By Corollary 4.3 and Proposition 6.2, the Coxeter graph $\Gamma_{VB,n}$ is of type $K(\pi,1)$. By \cite{ChaDav1}, Corollary 1.4.2, it follows that $K_n$ has a classifying space of dimension $d$ (which is  $\Omega/W$), and the cohomological dimension of $K_n$ is $d$. It remains to show that $d = n-1$.

\bigskip\noindent
Let $X \in \SS^f$. We say that a given {\it $i \in \{1, \dots, n\}$ occurs} in $X$ if there is some $j \in \{1, \dots, n\}$, $j \neq i$, such that either $x_{i,j} \in X$, or $x_{j,i} \in X$. Observe that an $i \in \{1, \dots, n\}$ cannot occur in $X$ more than twice, and that there exists an $i \in \{1, \dots, n\}$ which either occurs once, or never occurs (otherwise all the vertices of $\Gamma_X$ would be of valence $2$). This implies that $|X| \le n-1$. On the other hand, $Y=\{x_{1,2}, x_{2,3}, \dots, x_{n-1,n}\} \in \SS^f$ and $|Y|=n-1$.

\bigskip\noindent
The third part of Corollary 6.3 is a straightforward consequence of the second one.
\qed

\bigskip\noindent
{\bf Remark.}
Bardakov announced in \cite{Barda1} a solution to the word problem for virtual braid groups. However, Lemma 6 in \cite{Barda1}, which is the key point in the algorithm (as well as in the paper) is false. On the other hand, another solution to the word problem for virtual braid groups will appear in the forthcoming paper \cite{GodPar1}.

\bigskip\noindent
{\bf Remark.}
A classifying space $QC_n$ for the pure virtual braid group $VP_n$ is constructed in \cite{BaEnEtRa1}. This space is a CW complex of dimension $n-1$, thus the cohomological dimension of $VP_n$ is $\le n-1$, and $VP_n$ is torsion free. On the other hand, $VP_n$ contains the Artin pure braid group $P_n$ whose cohomological dimension is known to be $n-1$, thus the cohomological dimension of $VP_n$ is exactly $n-1$. This gives an alternative proof to Corollary 6.3.3.



\clearpage

\bigskip\bigskip\noindent
{\bf Eddy  Godelle,}

\smallskip\noindent
Université de Caen, Laboratoire LMNO, UMR 6139 du CNRS, Campus II, 14032 Caen cedex, France. 

\smallskip\noindent
E-mail: {\tt eddy.godelle@unicaen.fr}

\bigskip\noindent
{\bf Luis Paris,}

\smallskip\noindent 
Université de Bourgogne, Institut de Mathématiques de Bourgogne, UMR 5584 du CNRS, B.P. 47870, 21078 Dijon cedex, France.

\smallskip\noindent
E-mail: {\tt lparis@u-bourgogne.fr}



\begin{thebibliography}{99}

\bibitem{AbrBro1}
{\bf P. Abramenko, K. S. Brown.}
{\it Buildings. Theory and applications.}
Graduate Texts in Mathematics, 248. Springer, New York, 2008.

\bibitem{Allco1}
{\bf D. Allcock.}
{\it Braid pictures for Artin groups.}
Trans. Amer. Math. Soc. {\bf 354} (2002), no. 9, 3455--3474. 

\bibitem{Altob1}
{\bf J. A. Altobelli.}
{\it The word problem for Artin groups of FC type.}
J. Pure Appl. Algebra  {\bf 129}  (1998),  no. 1, 1--22. 

\bibitem{AltCha1}
{\bf J. A. Altobelli, R. Charney.}
{\it A geometric rational form for Artin groups of FC type.}
Geom. Dedicata  {\bf 79}  (2000),  no. 3, 277--289.

\bibitem{Artin1}
{\bf E. Artin.} 
{\it Theorie der Zöpfe.}
Abhandlungen Hamburg  {\bf 4} (1925), 47-72.

\bibitem{Artin2}
{\bf E. Artin.}
{\it Theory of braids.}
Ann. of Math. (2)  {\bf 48} (1947), 101--126.

\bibitem{Barda1}
{\bf V. G. Bardakov.}
{\it The virtual and universal braids.}
Fund. Math.  {\bf 184}  (2004), 1--18.

\bibitem{BarBel1}
{\bf V. G. Bardakov, P. Bellingeri.}
{\it Combinatorial properties of virtual braids.}
Topology Appl.  {\bf 156} (2009),  no. 6, 1071--1082.

\bibitem{BaEnEtRa1}
{\bf L. Bartholdi, B. Enriquez, P. Etingof, E. Rains.}
{\it Groups and Lie algebras corresponding to the Yang-Baxter equations.}
J. Algebra  {\bf 305}  (2006),  no. 2, 742--764.

\bibitem{Bourb1}
{\bf N. Bourbaki.}
{\it Eléments de mathématique. Fasc. XXXIV. Groupes et algèbres de Lie. Chapitre IV: Groupes de Coxeter et systèmes de Tits. Chapitre V: Groupes engendrés par des réflexions. Chapitre VI: Systèmes de racines.}
Actualités Scientifiques et Industrielles, No. 1337, Hermann, Paris, 1968. 

\bibitem{Brids1}
{\bf M. R. Bridson.}
{\it Geodesics and curvature in metric simplicial complexes.}
Group theory from a geometrical viewpoint (Trieste, 1990),  373--463, World Sci. Publ., River Edge, NJ, 1991. 

\bibitem{BriHae1}
{\bf M. R. Bridson, A. Haefliger.}
{\it Metric spaces of non-positive curvature.}
Grundlehren der Mathematischen Wissenschaften, 319. Springer-Verlag, Berlin, 1999. 

\bibitem{Bries1}
{\bf E. Brieskorn.}
{\it Die Fundamentalgruppe des Raumes der regulären Orbits einer endlichen komplexen Spiegelungsgruppe.}
Invent. Math. {\bf 12} (1971), 57--61. 

\bibitem{Bries2}
{\bf E. Brieskorn.}
{\it Sur les groupes de tresses [d'après V. I. Arnol'd].} 
Séminaire Bourbaki, 24ème année (1971/1972), Exp. No. 401, pp. 21--44. Lecture Notes in Math., Vol. 317, Springer, Berlin, 1973.

\bibitem{BriSai1}
{\bf E. Brieskorn, K. Saito.}
{\it Artin-Gruppen und Coxeter-Gruppen.}
Invent. Math. {\bf 17} (1972), 245--271.

\bibitem{Brown1}
{\bf K. S. Brown.}
{\it Buildings.}
Springer-Verlag, New York, 1989. 

\bibitem{Charn1}
{\bf R. Charney.}
{\it Artin groups of finite type are biautomatic.}
Math. Ann.  {\bf 292}  (1992),  no. 4, 671--683.

\bibitem{Charn2}
{\bf R. Charney.}
{\it Geodesic automation and growth functions for Artin groups of finite type.}
Math. Ann.  {\bf 301} (1995),  no. 2, 307--324.

\bibitem{ChaDav1}
{\bf R. Charney, M. W. Davis.}
{\it The $K(\pi,1)$-problem for hyperplane complements associated to infinite reflection groups.}
J. Amer. Math. Soc. {\bf 8} (1995), no. 3, 597--627. 

\bibitem{ChaDav2}
{\bf R. Charney, M. W. Davis.}
{\it Finite $K(\pi, 1)$s for Artin groups.}
Prospects in topology (Princeton, NJ, 1994),  110--124, Ann. of Math. Stud., 138, Princeton Univ. Press, Princeton, NJ, 1995. 

\bibitem{ChaPei1}
{\bf R. Charney, D. Peifer.}
{\it The $K(\pi,1)$-conjecture for the affine braid groups.}
Comment. Math. Helv.  {\bf 78}  (2003),  no. 3, 584--600.

\bibitem{Davis1}
{\bf M. W. Davis.}
{\it The geometry and topology of Coxeter groups.}
London Mathematical Society Monographs Series, 32. Princeton University Press, Princeton, NJ, 2008.

\bibitem{Delig1}
{\bf P. Deligne.}
{\it Les immeubles des groupes de tresses généralisés.} 
Invent. Math. {\bf 17} (1972), 273--302. 

\bibitem{Dieck1}
{\bf T. tom Dieck.}
{\it Categories of rooted cylinder ribbons and their representations.}
J. Reine Angew. Math. {\bf 494} (1998), 35--63.
 
\bibitem{Digne1}
{\bf F. Digne.}
{\it Présentations duales des groupes de tresses de type affine $\widetilde A$.}
Comment. Math. Helv.  {\bf 81}  (2006),  no. 1, 23--47. 

\bibitem{EllSko1}
{\bf G. Ellis, E. Sköldberg.}
{\it The $K(\pi,1)$ conjecture for a class of Artin groups.}
Comment. Math. Helv. {\bf 85} (2010), no. 2, 409--415.

\bibitem{FadNeu1}
{\bf E. Fadell, L. Neuwirth.}
{\it Configuration spaces.}
Math. Scand.  {\bf 10} (1962), 111--118.

\bibitem{Garsi1}
{\bf F.A. Garside.}
{\it The braid group and other groups.}
Quart. J. Math. Oxford Ser. (2)  {\bf 20} (1969), 235--254. 

\bibitem{GodPar1}
{\bf E. Godelle, L. Paris.}
{\it Locally Garside groups of FC type and the virtual braid group.}
In preparation.

\bibitem{Gromo1}
{\bf M. Gromov.}
{\it Hyperbolic groups.}
Essays in group theory, 75--263,
Math. Sci. Res. Inst. Publ., 8, Springer, New York, 1987. 

\bibitem{Kamad1}
{\bf S. Kamada.}
{\it Braid presentation of virtual knots and welded knots.}
Osaka J. Math.  {\bf 44}  (2007),  no. 2, 441--458. 

\bibitem{Kauff1}
{\bf L. H. Kauffman.}
{\it Virtual knot theory.}
European J. Combin.  {\bf 20}  (1999),  no. 7, 663--690.

\bibitem{KauLam1}
{\bf L. H. Kauffman, S. Lambropoulou.}
{\it Virtual braids and the $L$-move.}
J. Knot Theory Ramifications  {\bf 15}  (2006),  no. 6, 773--811.

\bibitem{KenPei1}
{\bf R. P. Kent, D. Peifer.}
{\it A geometric and algebraic description of annular braid groups.}
Internat. J. Algebra Comput. {\bf 12} (2002), no. 1-2, 85--97. 

\bibitem{Lek1}
{\bf H. van der Lek.}
{\it The homotopy type of complex hyperplane complements.}
Ph. D. Thesis, Nijmegen, 1983.

\bibitem{NibRee1}
{\bf G. A. Niblo, L. D. Reeves.}
{\it The geometry of cube complexes and the complexity of their fundamental groups.}
Topology  {\bf 37} (1998), no. 3, 621--633.

\bibitem{Okone1}
{\bf C. Okonek.}
{\it Das $K(\pi ,\,1)$-Problem für die affinen Wurzelsysteme vom Typ $A_{n}$, $C_{n}$.}
Math. Z. {\bf 168} (1979), no. 2, 143--148. 

\bibitem{Raben1}
{\bf L. Rabenda.}
Mémoire de DEA (Master thesis), Université de Bourgogne, 2003.

\bibitem{Salve1}
{\bf M. Salvetti.}
{\it The homotopy type of Artin groups.}
Math. Res. Lett.  {\bf 1} (1994),  no. 5, 565--577. 


\bibitem{Tits1}
{\bf J. Tits.}
{\it Le problème des mots dans les groupes de Coxeter.}
Symposia Mathematica (INDAM, Rome, 1967/68), Vol. 1, pp. 175--185, Academic Press, London, 1969. 


\bibitem{Weil1}
{\bf A. Weil.}
{\it Sur les théorèmes de de Rham.} 
Comment. Math. Helv. {\bf 26} (1952), 119--145. 

\end{thebibliography}
\end{document}